\documentclass[a4paper,10pt]{article}
\usepackage{authblk}
\usepackage[pdftex]{graphicx}
\usepackage{subcaption}

\usepackage[ngerman,english]{babel}
\usepackage[latin1]{inputenc}
\usepackage{csquotes}
\usepackage{amsfonts,amsmath,amsthm}
\usepackage[titletoc,title]{appendix}
\usepackage[backend=bibtex8,doi=false,eprint=false,firstinits=true,isbn=false,style=numeric-comp,url=false,maxnames=99]{biblatex}
\bibliography{bibl.bib}
\DeclareRedundantLanguages{english,german,french}{english,german,ngerman,french}

\usepackage{cases}
\usepackage{mathabx}
\usepackage{xfrac}
\usepackage{enumerate}
\usepackage{fancyhdr}
\usepackage{color}
\usepackage[colorlinks]{hyperref}
\definecolor{blue75}{rgb}{0,0,.75}
\definecolor{green75}{rgb}{0,.75,0}
\hypersetup{colorlinks=true, urlcolor=blue75,linkcolor=blue75,citecolor=green75,pdfstartview=FitB,bookmarksopen=true,bookmarksopenlevel=1}
\usepackage[a4paper, left=2.5cm, right=2.5cm, top=2.5cm,bottom=2cm]{geometry}
\usepackage{constants}
\newcommand{\parenthezises}[1]{\arabic{#1}}
\newconstantfamily{C}{
symbol=C,
format=\parenthezises,
reset={section}
}
\newconstantfamily{M}{
symbol=M,
format=\parenthezises,
reset={section}
}
\usepackage{enumitem}

\usepackage{graphicx}
\graphicspath{ {images/} }
\usepackage{wrapfig}
\usepackage{figbib}
\allowdisplaybreaks
\begin{document}
%%%Numbers%%%
\newcommand{\R}{\mathbb{R}}
\newcommand{\N}{\mathbb{N}}
%%%Operators%%%
\def\diam{\operatorname{diam}}
\def\dist{\operatorname{dist}}
\def\ess{\operatorname{ess}}
\def\inner{\operatorname{int}}
\def\osc{\operatorname{osc}}
\def\sign{\operatorname{sign}}
\def\supp{\operatorname{supp}}
%%%Spaces%%%
\newcommand{\BMO}{BMO\left(\Omega\right)}
\newcommand{\LOne}{L^{1}\left(\Omega\right)}
\newcommand{\LTwo}{L^{2}\left(\Omega\right)}
\newcommand{\Lq}{L^{q}\left(\Omega\right)}
\newcommand{\Lp}{L^{p}\left(\Omega\right)}
\newcommand{\LInf}{L^{\infty}\left(\Omega\right)}
\newcommand{\HOneO}{H^{1,0}\left(\Omega\right)}
\newcommand{\HTwoO}{H^{2,0}\left(\Omega\right)}
\newcommand{\HOne}{H^{1}\left(\Omega\right)}
\newcommand{\HTwo}{H^{2}\left(\Omega\right)}
\newcommand{\HmOne}{H^{-1}\left(\Omega\right)}
\newcommand{\HmTwo}{H^{-2}\left(\Omega\right)}

\newcommand{\LlogL}{L\log L(\Omega)}

\newcommand{\abs}{\\[0.5ex]}
\newcommand{\cb}{\color{blue}}
\newcommand{\veps}{v_{\varepsilon}}
\newcommand{\ceps}{c_{\varepsilon}}
\newcommand{\meps}{m_{\varepsilon}}

% \newcommand{\qn}{H^2_{loc}\left(\R^+_0,H^{-\sfrac{1}{2}}\left(\Gamma\right)\right)}
%%%Theorems%%%
\newtheorem{Theorem}{Theorem}[section]
\newtheorem{Assumption}[Theorem]{Assumptions}
\newtheorem{Corollary}[Theorem]{borollary}
\newtheorem{Convention}[Theorem]{bonvention}
\newtheorem{Definition}[Theorem]{Definition}
\newtheorem{Example}[Theorem]{Example}
\newtheorem{Lemma}[Theorem]{Lemma}
\newtheorem{Notation}[Theorem]{Notation}
\newtheorem{Remark}[Theorem]{Remark}
%%%Settings%%%
\numberwithin{equation}{section}

\title{Global existence for a degenerate  haptotaxis model of tumor invasion under the go-or-grow dichotomy hypothesis}
\author{Anna~Zhigun, Christina~Surulescu, and Alexander Hunt}
\renewcommand\Affilfont{\itshape\small}
\affil{Technische Universität Kaiserslautern, Felix-Klein-Zentrum für Mathematik\\ Paul-Ehrlich-Str. 31, 67663 Kaiserslautern, Germany\\
  e-mail: {\{zhigun,surulescu,hunt\}@mathematik.uni-kl.de}}
\date{}
\maketitle
\begin{abstract}%\noindent
We propose and study a strongly coupled PDE-ODE-ODE system modeling cancer cell invasion through a tissue network 
under the go-or-grow hypothesis asserting that cancer cells can either move or proliferate. Hence our setting features 
two interacting cell populations with their mutual transitions and involves tissue-dependent degenerate diffusion and 
haptotaxis for the moving subpopulation. The proliferating cells and the tissue evolution are characterized by way of ODEs 
for the respective densities. We prove the global existence of weak solutions and  illustrate the model behaviour by 
numerical simulations in a two-dimensional setting.
MSC 2010:
% 35B65, %Smoothness and regularity of solutions
35B45, %A priori estimates 
35D30, %Weak solutions
% 35D35, %Strong solutions
35K20, %Initial-boundary value problems for second-order parabolic equations
35K51, %Initial-boundary value problems for second-order parabolic systems 
% 35K55, %Nonlinear parabolic equations
35K57, %Reaction-diffusion equations
35K59, %Quasilinear parabolic equations
35K65, %Degenerate parabolic equations
35Q92,  %PDEs in connection with biology and other natural sciences
92C17. %Cell movement (chemotaxis, etc.)
\end{abstract}

\section{Introduction}\label{intro}
One of the essential characteristics of a tumor is its heterogeneity. The cells forming the neoplastic tissue often have different phenotypes, morphologies, and functions, and can switch between these in response to intra- and/or extracellular influences like e.g., genetic change, acidity of the peritumoral region, availability of nutrients and/or space, applied therapeutic agents etc., see e.g. \cite{kleppe-levine,marusyk-polyak,michor-weaver}. Tumor heterogeneity is tightly connected to compromised treatment response \cite{fidler,heppner} and is  already manifested at the migrating stage of tumor development. Indeed, one of the main features of tumor development and invasion is the ability of cancer cells to migrate and spread into the normal tissue, whereby they experience different migratory phenotypes (e.g., amoeboid vs. mesenchymal). Furthermore, experimental evidence revealed that several types of tumor cells (including glioma, breast cancer cells, and melanoma) defer their proliferation while migrating and vice versa \cite{giese-etal-96,giese-etal03,jerby-etal,widmer}, corresponding to the so-called go-or-grow dichotomy. The differentiated response of tumor cells to treatment is a main cause of radio- and chemotherapeutical failure; indeed, it is largely accepted that cells with a highly proliferating phenotype are more sensitive to therapy, whereas the migratory phenotype is attended by reduced treatment sensitivity, see e.g., \cite{lefrank,moore,steinbach-weller} and the references therein. \abs
\noindent
Motivated by the above mentioned facts we propose in this paper a model for tumor cell invasion in which we account for the go-or-grow hypothesis and distinguish between migrating and proliferating (hence non-moving) cells. Several continuum mathematical models relying on the go-or-grow behavior of tumor cells and explicitly accounting for the two subpopulations of migrating and proliferating cells, respectively, have been  considered e.g., in \cite{gerlee-nelander,saut} and featured reaction-(cross-)diffusion(-chemotaxis) equations. Using a two-component continuous-time random walk along with a probabilistic approach based thereupon and involving switching with exponentially distributed waiting times between the proliferation and migration phenotypes, Fedotov \& Iomin deduced in \cite{fedotov-iomin} an ODE-PDE system for the macroscopic dynamics of the two types of cancer cell densities, supporting the idea of tumor cells subdiffusivity instead of the more common Fickian diffusivity. In \cite{chauviere-preziosi-byrne} Chauviere et al. used a mesoscopic description of the two cell subpopulations to deduce by an appropriate scaling a system of two coupled reaction-diffusion equations for their  macrolevel behavior. Still in that context, starting from mesoscopic equations for the two cell subpopulations and coupling them with subcellular level dynamics in \cite{eks15,hs15} the authors obtained by parabolic scalings macroscopic equations characterizing the evolution of the overall tumor burden for a glioma invasion model. The resulting equations carried in their coefficients the information from the lower modeling scales (both subcellular and mesoscopic) and allowed DTI-based predictions about the tumor extent and simulation-based therapy outcomes. The haptotaxis term obtained in those macroscopic equations was a direct consequence of accounting for the subcellular receptor binding dynamics in the mesoscale evolution of the cancer cell densities. By using the equlibrium of fluxes and some ideas from \cite{mss-15,MSS_Acid}, in \cite{SSU} was introduced a multiscale model for macroscopic tumor invasion and development complying to the go-or-grow dichotomy and including subcellular dynamics of receptor binding to fibers of the underlying extracellular matrix (ECM). Our  model in this paper extends in a certain way the previous setting in \cite{SSU} by allowing the diffusion coefficient to degenerate and by paying increased attention to the haptotactic sensitivity function; however, neither therapy effects nor multiscality issues are addressed here.\abs
\noindent
While there is a vast literature concerning the mathematical analysis of reaction-diffusion-taxis equations, problems with degenerate diffusion and taxis have been less  investigated. However, during the last decade more such references became available; they 
describe the dynamics of a cell population in response to a chemoattractant \cite{eberl-efendiev-zhigun-14,lw-05,wang-winkler-wrzosek-12}, moving up the 
gradient of an insoluble signal (haptotaxis) \cite{ZSU,winkler-surulescu}, or performing both chemo- and haptotaxis 
\cite{li-lankeit,tao-winkler-11,zheng-etal-16}.
%, the latter with constant tactic sensitivities.  
Thereby, the type of degeneracy is a particularly relevant feature for the difficulty of the problem, especially for 
systems coupling ODEs with PDEs, as is the case when considering haptotaxis. In \cite{li-lankeit,tao-winkler-11,zheng-etal-16} 
the diffusion coefficients depend nonlinearly on the solution and the tactic sensitivities are constants. For these problems the 
global well posedness was obtained, along with boundedness properties of the solutions. The model proposed in \cite{ZSU} 
involves a diffusion coefficient which can degenerate due to each of the solution components (density of cells and of 
ECM fibers, respectively): moreover, the haptotactic sensitivity is a nonlinear function of the ECM density. The 
1D model in \cite{winkler-surulescu} was motivated by the deduction of macroscopic equations from a mesoscopic setting for 
brain tumor invasion also accounting for subcellular dynamics; it features a reaction-diffusion-transport-haptotaxis equation 
for the tumor cell density coupled with an ODE for the density of tissue fibers. The strong degeneracy of the diffusion and 
haptotaxis coefficients is attained by way of a function only depending on the position and not on the solution itself. 
Whereas the global existence of weak solutions was shown for these models, the boundedness and uniqueness issues remain open. 
The same applies to the mathematical setting considered in this work and presented in detail in the following {\it 
Section \ref{model}}. The rest of the paper is organized as follows: {\it Section \ref{not}} introduces some basic 
notations, {\it Section \ref{problem}} settles the problem and states the main result consisting in the 
global existence of a weak solution to the system in {\it Section \ref{model}}, to be followed by several steps 
towards its proof. Thus, {\it Section \ref{aproxi}} introduces a sequence of non-degenerate approximations of 
the actual problem and {\it Section \ref{apriori}} is concerned with deducing some apriori estimates to be used in 
{\it Section \ref{existence}} for the convergences necessary to prove the result announced in {\it Section \ref{problem}}. 
Finally, in \textit{Section \ref{numerics}} we perform some numerical simulations in order to illustrate the model behavior 
and we also comment on the obtained results.

\section{The model}\label{model}
% In our previous work \cite{ZSU}, we introduced the following system:
% \begin{subequations}\label{hapto}
% \begin{alignat}{3}
%  &\partial_t c=\nabla\cdot\left(\frac{\kappa_m vc}{1+vc}\nabla c-\frac{\kappa_v c}{(1+v)^2}\nabla v\right)+\mu_cc(1-c-\eta v)&&\text{ in }\R^+\times\Omega,\label{haptoc}\\
%  &\partial_t v=\mu_vv(1-v)-\lambda vc&&\text{ in }\R^+\times\Omega,\label{haptov}\\
%  &\frac{\kappa_m vc}{1+vc}\partial_{\nu} c-\frac{\kappa_v c}{(1+v)^2}\partial_{\nu} v=0&&\text{ in }\R^+\times\partial\Omega,\label{haptobc}\\
%  &c(0)=c_0,\ v(0)=v_0 &&\text{ in }\Omega,
% \end{alignat}
% \end{subequations}
% where $\eta,\lambda,\kappa_m, \kappa_v, \mu_c,\mu_v$ are some positive constants. System \eqref{hapto} can serve as a model prototype for the situation when just two components are present: the cancer cells  density 
% $c$ and the density $v$ of ECM tissue fibers, both depending on time and position on a smooth bounded domain $\Omega \subset \R^N$. In this work, we assume that the total cells density $c$ can be decomposed in two components: the density of the moving cells $m$ and that of proliferating cells $p$, so that
% \begin{align}
%  c=m+p.\nonumber%\label{cmcp}
% \end{align}
Based on the models in \cite{SSU,ZSU} we introduce here a PDE-ODE-ODE system characterizing the macroscopic dynamics of a tumor in interaction with the 
surrounding tissue  in accordance with the go-or-grow dichotomy.  The latter means that the tumor is assumed to be made up of two types of cells, which are either moving or mitotic and non-motile, whereby mutual transitions between the two phenotypes take place. Our model thus reads:
\begin{subequations}\label{haptoGG}
\begin{alignat}{3}
 &\partial_t m=-\alpha m+\beta v p+\nabla\cdot\left(\frac{\kappa_m vc}{1+vc}\nabla m
 -\frac{\kappa_v m}{(1+v)^2}\nabla v\right)&&\text{ in }\R^+\times\Omega,\label{cm}\\
 &\partial_t p=\alpha m-\beta vp+\mu_pp(1-c-\eta v)&&\text{ in }\R^+\times\Omega,\label{cp}\\
 &\partial_t v=\mu_vv(1-v)-\lambda vm&&\text{ in }\R^+\times\Omega,\label{v}\\
 &\frac{\kappa_m vc}{1+vc}\partial_{\nu} m-\frac{\kappa_v m}{(1+v)^2}\partial_{\nu} v=0&&\text{ in }\R^+\times\partial\Omega,\label{bc}\\
 &m(0)=m_{0},\ p(0)=p_{0},\ v(0)=v_0 &&\text{ in }\Omega,
\end{alignat}
\end{subequations}
where $m$ and $p$ denote the densities of moving and proliferating cells, respectively, $v$ is the density of ECM fibers, all depending on time and position on a smooth bounded domain $\Omega \subset \R^N$. The positive constants $\alpha,\beta $ denote the transition rates between the two subpopulations, $\eta >0$ is a constant 
scaling the concurrence with normal tissue in the proliferation process, $\kappa_m, \kappa_v$ are positive constants scaling the diffusion and the haptotactic sensitivity, $\lambda >0$ is the decay rate of ECM due to interactions with (mesenchymally) motile cells, and $\mu_p, \mu_v$ are growth rates for the tumor cells and the tissue, respectively. The total tumor burden is assessed by
$$c(t,x):=m(t,x)+p(t,x).$$
Thus, system \eqref{haptoGG} includes a degenerate parabolic PDE for the moving and an ODE for the proliferating tumor cells, together with an ODE for the tissue density, supplemented by the initial and the 'no-flux' boundary conditions. The latter complies with the fact that cancer cells do not leave the tissue hosting the original tumor. As in \cite{ZSU}, the diffusion coefficient in the equation for moving cells is nonlinear and can degenerate due to either tissue or  tumor cell densities. The haptotaxis coefficient is nonlinear as well; its form is motivated by the microlocal cell-tissue interactions (as explained in \cite{ZSU}) and whence keeps a flavor of multiscality, also in a rather indirect fashion, as our system  \eqref{haptoGG} is purely macroscopic. For explicit multiscale effects we refer to the related model in \cite{SSU}. As observed there, the analysis done for a system involving a single population of cancer cells (hence without accounting for tumor heterogeneity in the sense mentioned above) does not directly carry over to a model discerning between moving and proliferating cells. One of the difficulties comes from the switching between the two populations, as the moving cells act on the one side as source for the proliferating ones, and on the other side as decay term for themselves and for the tissue. Another complicacy is due to the supplementary ODE for the proliferating cells, which -like the equation for tissue dynamics- lacks  space derivatives, which was already a challenge in the more classical haptotaxis settings. Here the degenerate diffusion renders the problem even more complex.\abs
\noindent

\section{Basic notation and functional spaces}\label{not}
We denote the Lebesgue measure of a set $A$ by $|A|$ and by $\inner A$ its interior.\\
\noindent
Partial derivatives, in both classical and distributional sense, with respect to variables  $t$ and $x_i$, will be denoted respectively by $\partial_t$ 
and $\partial_{x_i}$. Further, $\nabla$, $\nabla\cdot$ and $\Delta$ stand for the spatial gradient, divergence and Laplace operators, respectively. $\partial_{\nu}$ is 
the derivative with respect to the outward unit normal of $\partial\Omega$. 
% {\cb M\"ussen noch sagen, was $u_0$ ist (wir haben das auch im anderen Paper vergessen)}\abs
% \noindent
We assume the reader to be familiar with the standard Lebesgue and Sobolev spaces and their usual properties, as well as with the  more general $L^p$ spaces of functions with values in general Banach spaces and with anisotropic Sobolev spaces. In particular, we need the Banach space
\begin{align}
 &W^{-1,1}(\Omega):=\left\{u\in D'(\Omega)\left|\ u=u_0+\sum_{k=1}^N\partial_{x_i}u_i\text{ for some }u_i\in\LOne,\ i={0},\dots,N\right.\right\},\nonumber\\
 &||u||_{W^{-1,1}(\Omega)}:=\inf\left\{\sum_{k=0}^N||u_i||_1\left|\ u=u_0+\sum_{k=1}^N\partial_{x_i}u_i,\   u_i\in\LOne,\ i={0},\dots,N\right.\right\}.\nonumber
\end{align} 
We will also make use of the %'logarithmic' Lebesgue space
Zygmund space \cite[Chapter 6, Definition 6.1]{BeS88}
% \begin{align}
%  \LlogL:=\left\{u\in\LOne|\ \int_{\Omega}\chi_{\{|u|>1\}}|u|\log|u|\,dx<\infty\right\}.\label{LlogL}
% \end{align}
\begin{align}
 \LlogL:=\left\{u\in\LOne\left|\ \ \int_{\Omega}M(u)\,dx<\infty\right\},\text{ where }M(u):=\chi_{\{|u|>1\}}|u|\log|u|.\right.\nonumber%\label{LlogL}
\end{align}
For $p\in[1,\infty]\backslash\left\{2\right\}$, we write $||\cdot||_p$ in place of the $||\cdot||_{\Lp}$-norm. Throughout the paper, $||\cdot||$ and ($(u,v)$) denote  the standard $\LTwo$-norm  and scalar product, respectively.% and  for $\int_{\Omega}u(x)v(x)\,dx$.
% , while $\left<\cdot,\cdot\right>$ is reserved for the duality pairing between 
% $W^{1,\infty}(\Omega)$ and its dual $(W^{1,\infty}(\Omega))'$. 
\\
\noindent
Finally, we make the following useful convention: For all indices $i$, the quantity $C_i$ denotes a non-negative constant or, alternatively, a non-negative function, 
which is non-decreasing in each of its arguments.

\section{Problem setting and main result}\label{problem}
In this section we propose a definition of weak solutions to system \eqref{haptoGG} and state our main result under the following assumptions:
\begin{Assumption}[Initial data]\label{ini}~%Let 
\begin{enumerate} 
 \item 
  $m_{0}\geq0,\ m_{0}\notequiv 0, \ m_{0}\in\LlogL$;
  \item $p_{0}\geq0,\ p_{0}\notequiv 0, \  p_{0}\in\LInf$;
  \item
  $0\leq v_0\leq1,\ v_0\notequiv 0,1,\ v_0^{\frac{1}{2}}\in \HOne$.
\end{enumerate}
\end{Assumption}

\noindent
 The major challenge of model \eqref{haptoGG} lies in the fact that the diffusion coefficient in equation \eqref{cm} degenerates 
%  at $m=0$ and, moreover, 
 at $v=0$. The latter seems to make it impossible to obtain an a priori estimate for the gradient of  $\varphi(m)$ in some 
 Lebesgue space for any  smooth, strictly increasing function  $\varphi$.  As a workaround, we are forced to consider an auxiliary function 
%  $\ln\left(1+v^{\frac{1}{2}}c\right)$ 
 which involves {\it both} $m$ and $v$ and whose gradient we are able to estimate.\\
% \begin{align}
%   \ln\left(1+v^{\frac{1}{2}}c\right).\nonumber
% \end{align}
\noindent
This leads us to the following definition of weak solutions to  \eqref{haptoGG}: 
\begin{Definition}[Weak solution]\label{Defweak}
Let $m_0,p_0,v_0$ satisfy {\it Assumptions~\ref{ini}}.
We call a triple of functions $m,p:\R^+_0\times\overline{\Omega}\rightarrow\R^+_0$, $v:\R^+_0\times\overline{\Omega}\rightarrow[0,1]$ a global weak solution of
 \eqref{haptoGG} if for all $0<T<\infty$ it holds that
\begin{enumerate}
  \item $m\in L^{\infty}(0,T;\LOne)$;%, $\partial_t m\in L^1(0,T;(W^{1,\infty}(\Omega))')$;
  %$\partial_t c\in L^1(0,T;W^{-1,1}(\Omega)\cap (W^{1,\infty}(\Omega))')$;
  \item $p\in L^{\infty}(0,T;\LInf)$, $\partial_t p\in L^{1}(0,T;\LOne)$;
  \item $v^{\frac{1}{2}}\in L^{\infty}(0,T;\HOne)$, $\partial_t v\in L^{1}(0,T;\LOne)$;
%   \item $\ln\left(1+v^{\frac{1}{2}}m\right)\in L^{1}(0,T;W^{1,1}(\Omega))$, $\frac{ v^{\frac{1}{2}}c}{1+vc}\left(\left(1+v^{\frac{1}{2}}m\right)\nabla \ln\left(1+v^{\frac{1}{2}}m\right)-m\nabla v^{\frac{1}{2}}\right)\in L^1(0,T;\LOne)$;
  \item $\nabla\left(v^{\frac{1}{2}}(m+1)^{\frac{1}{2}}\right)$, $\left(\frac{ \kappa_m c}{1+vc}+\frac{\kappa_v}{1+v}\right)v^{\frac{1}{2}}(m+1)^{\frac{1}{2}}\left(\nabla\left(v^{\frac{1}{2}}(m+1)^{\frac{1}{2}}\right)-(m+1)^{\frac{1}{2}}\nabla v^{\frac{1}{2}}\right)\in L^1(0,T;\LOne)$, $\nabla\left(\int_0^t\frac{v}{1+v}m\,d\tau\right)\in L^{\infty}(0,T;\LTwo)$;
  \item $(m,p,v)$ satisfies equation \eqref{cm} and the boundary condition \eqref{bc} in the following weak sense: 
  \begin{align}
   &\int_{\Omega}m_0\varphi\,dx\psi(0)-\int_0^T\int_{\Omega}m\varphi\,dx\psi'\,dt\nonumber\\
   =&-\int_0^T\int_{\Omega}\left(\frac{ \kappa_m c}{1+vc}+\frac{\kappa_v}{1+v}\right)2v^{\frac{1}{2}}(m+1)^{\frac{1}{2}}\left(\nabla\left(v^{\frac{1}{2}}(m+1)^{\frac{1}{2}}\right)-(m+1)^{\frac{1}{2}}\nabla v^{\frac{1}{2}}\right)\cdot\nabla\varphi\psi\nonumber\\
   &\quad\quad\quad+\kappa_v\nabla\left(\int_0^t\frac{v}{1+v}m\,d\tau\right)\cdot\nabla\varphi\psi'+(-\alpha m+\beta v p)\varphi\psi\,dx\,dt\label{weakf}
  \end{align}
  for all $\varphi\in W^{1,\infty}(\Omega)$ and $\psi\in W^{1,\infty}(0,T)$ such that $\psi(T)=0$;
  \item $(m,p,v)$ satisfies equations \eqref{cp}-\eqref{v}  in $L^{1}(0,T;\LOne)$;
  \item % $m(0)=m_{0}$, 
  $p(0)=p_{0}$, $v(0)=v_0$.
 \end{enumerate}
\end{Definition}
{
\begin{Remark}[Weak formulation]\label{remweak}
By using the chain and product rules and (where necessary) partial integration over $\Omega$ and over $[0,t]$, it can be readily checked that  \eqref{weakf} is, indeed, a weak reformulation of \eqref{cm} and \eqref{bc}. Its somewhat nonstandard form is due to the fact that $\nabla m$ in the diffusion term and the taxis flux term  $\frac{\kappa_v m}{(1+v)^2}\nabla v$ might not exist even in $L^1_{loc}$-sense.  
\end{Remark}
}
\begin{Remark}[Initial conditions]
Since we are looking for solutions with 
\begin{align*}
% &m\in W^{1,1}((0,T);W^{-1,1}(\Omega)),\\
&p\in W^{1,1}(0,T;\LOne),\\
&v^{\frac{1}{2}}\in H^1(0,T;\LOne),                                                 \end{align*}
we have 
\begin{align*}
% &m\in C([0,T];W^{-1,1}(\Omega)),\\
&p\in C([0,T];\LOne),\\
&v^{\frac{1}{2}}\in C([0,T];\LOne).                                                 \end{align*}
Therefore, the initial conditions 7. in {\it Definition \ref{Defweak}} do make sense.
\end{Remark}
\noindent
Our main result reads:
\begin{Theorem}[Global existence]\label{maintheo}
Let $\Omega\subset\R^N$, $N\in\N$, be a smooth bounded domain and let $\alpha,\beta$, $\eta, \kappa_m, \kappa_v, \lambda, \mu_p, \mu_v$ be positive constants. 
 Then, for each triple of functions $(m_0,p_0,v_0)$ satisfying {\it Assumptions~\ref{ini}}
 there exists a global weak solution $(m,p,v)$ (in terms of {\it Definition~\ref{Defweak}})  to the system \eqref{haptoGG}.
\end{Theorem}
\noindent
The proof of {\it Theorem~\ref{maintheo}} is based on a suitable approximation of the degenerate system \eqref{haptoGG} by a family of systems with nondegenerate diffusion of the migrating cells,  
derivation of a set of a priori estimates which ensure necessary compactness and, finally, the passage to the limit. 
While the overall structure of the proof  
is a standard one for a haptotaxis system,  we encounter considerable difficulties in each of the three steps due to the previously mentioned degenerate diffusion 
in equation \eqref{cm}, due to the ODEs \eqref{cp}-\eqref{v} having  no diffusion at all (i.e., everywhere degenerate), and, finally, due to  the strong couplings. 

% In order to gain a better understanding of the underlying complexity, let us compare \eqref{haptoGG} to the following problem which was introduced in \cite{Kawasaki1997}:
% \begin{subequations}\label{oldhapto}
% \begin{alignat}{3}
%  &\partial_t c=\nabla\cdot\left(\frac{\kappa}{1+vc}\nabla c-\frac{\kappa vc}{(1+v)^2}\nabla v\right)+-\alpha m+\beta v p&&\text{ in }\R^+\times\Omega,\label{oldc}\\
%  &\partial_t v=\mu_vv(1-v)-\lambda vc&&\text{ in }\R^+\times\Omega,\label{oldv}
%  \end{alignat}
% \end{subequations}
%  and
% System \eqref{oldhapto} is a subsystem of the system studied in \cite{SSW}. 
% \begin{subequations}\label{nb_}
% \begin{alignat}{3}
% &\partial_t c=\nabla\cdot\left(\kappa vc\nabla c\right)+\theta f(v,c)&&\text{ in }\R^+\times\Omega,\label{b_}\\
% &\partial_t v=D_v\Delta v-f(v,c)&&\text{ in }\R^+\times\Omega.\label{n_}
% \end{alignat}
% \end{subequations}
% In \cite{HKWZ}, a generalisation of this model was studied analytically. It was observed there that 

\begin{Remark}[Notation for constants]
We make the following useful convention: 
The statement that a constant depends on the parameters of the problem means that it depends on the constants $\kappa, \mu_p, \eta, \mu_v$ and $\lambda$,
% and $\theta$ (see below),
the norms of the initial data $(m_0,p_0,v_0)$, the space dimension $N$, and the domain $\Omega$. This dependence on the parameters is subsequently {\bf not} indicated in an explicit way.
\end{Remark}

\section{Approximating problems}\label{aproxi}
In this section we introduce a family of non-degenerate approximations for problem \eqref{haptoGG}. 
For each relaxation parameter $\varepsilon=(\varepsilon_{1},\varepsilon_{2},\varepsilon_2)\in(0,1)^3$, the corresponding approximation reads  
\begin{subequations}\label{hapto2e}
\begin{alignat}{3}
  &\partial_t m_{\varepsilon}=-\alpha m_{\varepsilon}+\beta  v_{\varepsilon}p_{\varepsilon}+\varepsilon_1\Delta m_{\varepsilon}+\nabla\cdot\left(\frac{\kappa_m v_{\varepsilon}c_{\varepsilon}}{1+v_{\varepsilon}c_{\varepsilon}}\nabla m_{\varepsilon}-\frac{\kappa_v m_{\varepsilon}}{(1+v_{\varepsilon})^2}\nabla v_{\varepsilon}\right)%-\varepsilon_1 m_{\varepsilon}^{\theta}
  &&\text{ in }\R^+\times\Omega,\label{eq1e}\\
   &\partial_t p_{\varepsilon}=\alpha m_{\varepsilon}-\beta v_{\varepsilon} p_{\varepsilon}+\mu_pp_{\varepsilon}\left(1-c_{\varepsilon}-\eta v_{\varepsilon}\right)&&\text{ in }\R^+\times\Omega,\label{cpe}\\
 &\partial_t v_{\varepsilon}=\mu_vv_{\varepsilon}(1-v_{\varepsilon})-\lambda v_{\varepsilon}m_{\varepsilon}&&\text{ in }\R^+\times\Omega,\label{eq2e}\\
 &\varepsilon_1\partial_{\nu} m_{\varepsilon}+\frac{\kappa_m v_{\varepsilon}c_{\varepsilon}}{1+v_{\varepsilon}c_{\varepsilon}}\partial_{\nu} m_{\varepsilon}-\frac{\kappa_v m_{\varepsilon}}{(1+v_{\varepsilon})^2}\partial_{\nu} v_{\varepsilon}=0&&\text{ in }\R^+\times\partial\Omega,\label{bc1}\\
 &m_{\varepsilon}(0)=m_{\varepsilon_20},\ p_{\varepsilon}(0)=p_{\varepsilon_20},\ v_{\varepsilon}(0)=v_{\varepsilon_30} &&\text{ in }\Omega,\label{inie}
\end{alignat}
\end{subequations}
where 
\begin{align*}
 c_{\varepsilon}=m_{\varepsilon}+p_{\varepsilon}%\label{cmcpe}
\end{align*}
and the families $\{m_{\varepsilon_20}\}$, $\{p_{\varepsilon_20}\}$ and $\{v_{\varepsilon_30}\}$ of sufficiently smooth and nonnegative initial values are parameterized by $\varepsilon_2$ and $\varepsilon_3$, respectively.
They are yet to be %further
specified below in Subsection \ref{aid}.

\noindent
For each $\varepsilon\in(0,1)^3$, system  \eqref{hapto2e} has the form of a nondegenerate\footnote{in the sense that the parabolic PDE for the moving cells is nondegenerate} quasilinear haptotaxis system with respect to variables $m_{\varepsilon},p_{\varepsilon},v_{\varepsilon}$. { Thereby, the weak solutions can be defined similarly to {\it Definition~\ref{Defweak}}. In this case, 5. in Definition~\ref{Defweak} is replaced by
\begin{enumerate}
\item[5$^\prime$.]
$(m_{\varepsilon},p_{\varepsilon},v_{\varepsilon})$ satisfies equation \eqref{eq1e} and the boundary condition \eqref{bc1} in the following weak sense: 
  \begin{align}
   &\int_{\Omega}m_{\varepsilon_20}\varphi\,dx\psi(0)-\int_0^T\int_{\Omega}m_{\varepsilon}\varphi\,dx\psi'\,dt\nonumber\\
   =&\int_0^T\int_{\Omega}-\varepsilon\nabla m_{\varepsilon}\cdot\nabla\varphi\psi\,dxdt\nonumber\\
   &-\int_0^T\int_{\Omega}\left(\frac{ \kappa_m c_{\varepsilon}}{1+v_{\varepsilon}c_{\varepsilon}}+\frac{\kappa_v}{1+v_{\varepsilon}}\right)2v_{\varepsilon}^{\frac{1}{2}}(m_{\varepsilon}+1)^{\frac{1}{2}}\left(\nabla\left(v_{\varepsilon}^{\frac{1}{2}}(m_{\varepsilon}+1)^{\frac{1}{2}}\right)-(m_{\varepsilon}+1)^{\frac{1}{2}}\nabla v_{\varepsilon}^{\frac{1}{2}}\right)\cdot\nabla\varphi\psi\nonumber\\
   &\quad\quad\quad+\kappa_v\nabla\left(\int_0^t\frac{v_{\varepsilon}}{1+v_{\varepsilon}}m_{\varepsilon}\,d\tau\right)\cdot\nabla\varphi\psi'+(-\alpha m_{\varepsilon}+\beta v_{\varepsilon} p_{\varepsilon})\varphi\psi\,dx\,dt\label{weakfe}
  \end{align}
  for all $\varphi\in W^{1,\infty}(\Omega)$ and $\psi\in W^{1,\infty}(0,T)$ such that $\psi(T)=0$.
\end{enumerate}
}\noindent
The global existence of nonnegative weak solutions for system \eqref{hapto2e} can be obtained in a standard way.  
% {\cb Bitte schwache Formulierung explizit hinschreiben.}
We refer the reader to our proof in \cite{ZSU}  where we dealt with a similar situation. It is based on further regularizations, Amann's theory for abstract parabolic quasilinear systems \cite{Amann1}, and a priori estimates. We omit those details here. \\
It is clear that for $\varepsilon=0$ we regain - at least formally - the original degenerate haptotaxis system \eqref{haptoGG}. %{\cb Wir sollten sagen, was  $m_{\varepsilon_20}$ usw sind, damit klar wird wie man f\"ur $\varepsilon =0$ wieder $m_0$ usw. bekommt.} 
As it turns out (see the subsequent {\it Section~\ref{existence}}), a weak solution to \eqref{haptoGG} can be obtained as a limit of a sequence of solutions to \eqref{hapto2e}.\\

\noindent
In order to shorten the writing, we will sometimes use the following notation for the flux and reaction terms, respectively:
\begin{align}
 &q_{\varepsilon}:=\varepsilon_1\nabla m_{\varepsilon}+\frac{\kappa_m v_{\varepsilon}c_{\varepsilon}}{1+v_{\varepsilon}c_{\varepsilon}}\nabla m_{\varepsilon}-\frac{\kappa_v m_{\varepsilon}}{(1+v_{\varepsilon})^2}\nabla v_{\varepsilon},\label{qe}\\
 &f_{\varepsilon}:=-\alpha m_{\varepsilon}+\beta  v_{\varepsilon}p_{\varepsilon}\label{fe}
\end{align}
\subsection{Approximating initial data}\label{aid}
Our next step is to construct a suitable family of approximations to the initial data. Since we assume that $(m_0,p_0,v_0)$ satisfies {\it Assumptions~\ref{ini}}, there exists for each $\left(\varepsilon_2,\varepsilon_3\right)\in(0,1)^2$ an approximation triple  $\left(m_{\varepsilon_20},p_{\varepsilon_20},v_{\varepsilon_30}\right)$ with the following properties:
\begin{align}
 &m_{\varepsilon_20},p_{\varepsilon_20},v_{\varepsilon_30}^{\frac{1}{2}}\in W^{1,\infty}(\Omega),\label{e34bnd}\\
 &m_{\varepsilon_20},p_{\varepsilon_20}\geq0,\ 0\leq v_{\varepsilon_30}\leq1 \text{ in }\overline{\Omega},\ m_{\varepsilon_20},p_{\varepsilon_20}, v_{\varepsilon_30}\notequiv0,\label{mp0}\\
%  &||m_{\varepsilon_20}\ln m_{\varepsilon_20}||\leq2||m_{0}\ln m_{0}||_1,\\
 &\left\|M\left(m_{\varepsilon_20}\right)\right\|_1\leq2\left\|M(m_{0})\right\|_1,\\
 & \left\|\nabla v_{\varepsilon_30}^{\frac{1}{2}}\right\|\leq2\left\| v_0^{\frac{1}{2}}\right\|_{H^1(\Omega)},\\
 &||m_{\varepsilon_20}-m_{0}||_1\leq\varepsilon_2,\label{moe}\\
 &||p_{\varepsilon_20}-p_{0}||_{\infty}\leq\varepsilon_2,\label{poe}\\
 & \left\|v_{\varepsilon_30}^{\frac{1}{2}}-v_{0}^{\frac{1}{2}}\right\|\leq\varepsilon_3.\label{voe}
\end{align}
Recall that our aim is to pass to the limit for $\varepsilon \to 0$ in the approximating problem. Since equation \eqref{eq2e} is an ODE,  
the set $\{v(t,\cdot )=0\}$ is preserved in time (possibly up to some subsets of measure zero). Therefore, it turns out that we have to pay particular care at the set $\{v_{\varepsilon_30}=0\}$ whose interior should not shrink substantially with respect to $\{v_0=0\}$. Following the idea from \cite{ZSU}, we assume that
\begin{align}
 \left|\{v_0=0\}\backslash\inner\left\{v_{\varepsilon_30}=0\right\}\right|\leq\varepsilon_3.\label{star}
\end{align}
% where $\inner\left\{v_{\varepsilon_30}=0\right\}$ denotes the interior of the set $\left\{v_{\varepsilon_30}=0\right\}$. 
Indeed, to justify \eqref{voe} we recall here our argument from \cite{ZSU} for the convenience of the reader. 
Due to a Lusin property for Sobolev functions  \cite[Chapter 6, Theorem 6.14]{EvansGar}, there exists a function $\xi$ such that
\begin{align}
&\xi\in W^{1,\infty}(\Omega),\label{L1}\\
 & \left\|\xi\right\|_{H^1(\Omega)}\leq2\left\|v_0^{\frac{1}{2}}\right\|_{H^1(\Omega)},\label{L12}\\
 &\left|\left\{\xi\neq v_0^{\frac{1}{2}}\right\}\right|\leq\frac{\varepsilon_3}{4}.\label{L2}
\end{align}
We define  $$v_{\varepsilon_30}:=\left(\min\{\xi_+,1\}-\frac{\varepsilon_3}{2|\Omega|}\right)_+^2.$$ Let us check that $v_{\varepsilon_30}$ satisfies the above assumptions.  Indeed, due to \eqref{L1}-\eqref{L12}, we have that
\begin{align}
&v_{\varepsilon_30}^{\frac12}\in W^{1,\infty}(\Omega),\nonumber\\
 &\left\|\nabla v_{\varepsilon_30}^{\frac{1}{2}}\right\|\leq ||\nabla \xi||\leq2\left\|v_0^{\frac{1}{2}}\right\|_{H^1(\Omega)},\nonumber
 \end{align}
and
\begin{align}
 \left\|v_{\varepsilon_30}^{\frac{1}{2}}-v_{0}^{\frac{1}{2}}\right\| \leq &2\left|\left\{\xi\neq v_0^{\frac{1}{2}}\right\}\right|+\left\|\chi_{\left\{\xi= v_0^{\frac{1}{2}}\right\}}\left(\left(\xi-\frac{\varepsilon_3}{2|\Omega|}\right)_+-\xi\right)\right\|\nonumber\\
 \leq&\varepsilon_3.\nonumber
 \end{align}
Moreover, it holds that 
 \begin{align}
 &\{\xi=0\}\subset\left\{\min\{\xi_+,1\}<\frac{\varepsilon_3}{2|\Omega|}\right\}\subset\inner\left\{\min\{\xi_+,1\}\leq\frac{\varepsilon_3}{2|\Omega|}\right\}\cup\partial\Omega=\inner\{v_{\varepsilon_30}=0\}\cup\partial\Omega.\label{L4}
\end{align}
Combining \eqref{L2} and \eqref{L4}, we obtain  \eqref{star}.
\section{A priori estimates}\label{apriori}
In this section we establish, based on system \eqref{hapto2e}, several uniform a priori estimates for the functions  $m_{\varepsilon},p_{\varepsilon},v_{\varepsilon}$ and their combinations, which we will use in the existence proof (see {\it Section~\ref{existence}} below). Our calculations make use of the regularity which the solutions of \eqref{hapto2e} do have. While operating with the weak derivatives, we use the weak chain and product rules. Another way to justify the calculation is via further approximations, as was done in \cite{ZSU}. 
{ \subsection*{Uniform boundedness of $v_{\varepsilon}$}
Since the ODE \eqref{eq2e}  has the form
\begin{align}
 \partial_t v_{\varepsilon}=f_v(v_{\varepsilon},m_{\varepsilon})\nonumber
\end{align}
with $f_v(0,m)=0$, $f_v(1,m)\leq0$ for all $m\geq 0$, and the initial value satisfies $0\leq v_{\varepsilon_30}\leq1$ (compare \eqref{mp0}), we obtain using standard ODE theory that
\begin{align}
 0\leq v_{\varepsilon}\leq 1 \text{ in }(0,T)\times\Omega.\nonumber
\end{align}
holds a priori. Below we will use this simple estimate without referring to it explicitly.
}
\subsection*{Uniform boundedness of $p_{\varepsilon}$}
Equation \eqref{cpe} for $p_{\varepsilon}$ can be rewritten in the following way:
\begin{align}
 \partial_t p_{\varepsilon}=-\left(\mu_pp_{\varepsilon}-\alpha\right) m_{\varepsilon}-(\beta+\mu_p\eta) v_{\varepsilon}p_{\varepsilon}+\mu_pp_{\varepsilon}\left(1-p_{\varepsilon}\right).\label{cpe1}
\end{align}
Since $m_{\varepsilon},v_{\varepsilon}\geq0$, one readily obtains from \eqref{cpe1} using the Gronwall lemma  that
\begin{align}
 p_{\varepsilon}\leq C_p.\label{cpb}
\end{align}
\subsection*{Energy-type estimates}
We now turn to equation \eqref{eq2e} for $v_{\varepsilon}$. On both sides of \eqref{eq2e}, we divide by $v_{\varepsilon}^{\frac{1}{2}}(1+v_{\varepsilon})$ and then apply the gradient operator. Thus we obtain that
\begin{align}
\partial_t \nabla \int_0^{v_{\varepsilon}}\frac{1}{s^{\frac{1}{2}}(1+s)}\,ds=&%\varepsilon_{1}\Delta \nabla \psi(v_{\varepsilon})
-\lambda \frac{v_{\varepsilon}^{\frac{1}{2}}}{1+v_{\varepsilon}}\nabla m_{\varepsilon}-\frac{\lambda (1-v_{\varepsilon})m_{\varepsilon}+\mu_v(-1+4v_{\varepsilon}+v_{\varepsilon}^2)}{(1+v_{\varepsilon})^2}\nabla v_{\varepsilon}^{\frac{1}{2}}.
 \label{deq2e}
\end{align}
Further, we multiply \eqref{eq1e} by 
$\ln m_{\varepsilon}$ and \eqref{deq2e}  by $\frac{\kappa_v}{\lambda}\nabla \int_0^{v_{\varepsilon}}\frac{1}{s^{\frac{1}{2}}(1+s)}\,ds$ and integrate over $\Omega$ using partial integration and the boundary conditions where necessary. Adding the resulting identities together, we obtain after some calculation that
\begin{align}
 &\frac{d}{dt}\left(\left(1,m_{\varepsilon}\ln m_{\varepsilon}-m_{\varepsilon}\right)+\frac{2\kappa_v}{\lambda}\left(\frac{1}{(1+v_{\varepsilon})^2},\left|\nabla v_{\varepsilon}^{\frac{1}{2}}\right|^2\right)\right)%+\varepsilon_1\left(\frac{1}{\theta}\left(m_{\varepsilon}^{\theta},\ln m_{\varepsilon}^{\theta}\right)+\left\|\Delta \psi(v_{\varepsilon})\right\|^2\right)
 +\varepsilon_1\left\|\nabla m_{\varepsilon}^{\frac{1}{2}}\right\|^2+4\left(\frac{\kappa_mv_{\varepsilon}c_{\varepsilon}}{1+v_{\varepsilon}c_{\varepsilon}},\left|\nabla m_{\varepsilon}^{\frac{1}{2}}\right|^2\right)\nonumber\\
 &+\frac{2\kappa_v}{\lambda}\left(\lambda(1-v_{\varepsilon})m_{\varepsilon}+5\mu_vv_{\varepsilon}+\mu_vv_{\varepsilon}^2,\frac{\left|\nabla v_{\varepsilon}^{\frac{1}{2}}\right|^2}{(1+v_{\varepsilon})^3}\right)+\alpha (m_{\varepsilon},\ln m_{\varepsilon})\nonumber\\
 \leq&\beta(v_{\varepsilon}p_{\varepsilon},\ln m_{\varepsilon})+%\left(\mu_p c_{\varepsilon} \left(1-c_{\varepsilon}-\eta v_{\varepsilon}\right),\ln c_{\varepsilon} \right)
 \frac{2\mu_v\kappa_v}{\lambda}\left(\frac{1}{(1+v_{\varepsilon})^2},\left|\nabla v_{\varepsilon}^{\frac{1}{2}}\right|^2\right)\nonumber
%  \\
%  \leq&-\Cl{C4}\left(\chi_{\{c_{\varepsilon}>1\}},c_{\varepsilon}^2\ln c_{\varepsilon}\right)+\frac{2\mu_v\kappa}{\lambda}\left(\frac{1}{(1+v_{\varepsilon})^2},\left|\nabla v_{\varepsilon}^{\frac{1}{2}}\right|^2\right)+\Cl{C5}.
\end{align}
By using the Gronwall lemma and \eqref{cpb}, we thus arrive, for arbitrary $T\in\R^+$, at the estimates
\begin{align}
& \sup_{t\in[0,T]}\left(\chi_{\{m_{\varepsilon}>1\}},m_{\varepsilon}\ln m_{\varepsilon}\right)\leq\Cr{C2}(T),\label{a1}\\ & \sup_{t\in[0,T]}\left\|\nabla v_{\varepsilon}^{\frac{1}{2}}\right\|^2\leq\Cl{C2}(T),\label{a2}\\
&\int_0^T\left(\frac{v_{\varepsilon}c_{\varepsilon}}{1+v_{\varepsilon}c_{\varepsilon}},\left|\nabla m_{\varepsilon}^{\frac{1}{2}}\right|^2\right)\,dt\leq\Cr{C2}(T),\label{a4c}\\ 
& \int_0^T\left((1-v_{\varepsilon})m_{\varepsilon},\left|\nabla v_{\varepsilon}^{\frac{1}{2}}\right|^2\right)\,dt\leq\Cr{C2}(T),\label{a5}\\
 &\int_0^T\left(v_{\varepsilon}p_{\varepsilon},-\chi_{\{m_{\varepsilon}<1\}}\ln m_{\varepsilon}\right)\,dt\leq\Cr{C2}(T)\\
&\int_0^T\left\|\nabla m_{\varepsilon}^{\frac{1}{2}}\right\|^2\,dt\leq\varepsilon_{1}^{-1}\Cr{C2}(T).\label{a6}%\\ 
% & \int_0^T\left\|m_{\varepsilon}^{\theta}\ln m_{\varepsilon}^{\theta}\right\|_1\,dt\leq\varepsilon_{1}^{-1}\Cr{C2}(T),\label{a7}\\ 
% & \int_0^T\left\|\Delta \psi(v_{\varepsilon})\right\|^2\,dt\leq\varepsilon_{1}^{-1}\Cr{C2}(T).\label{a8}
\end{align}
Since $s\mapsto\frac{s}{1+s}$ is a monotonically increasing function, \eqref{a4c} yields that
\begin{align}
 \int_0^T\left(\frac{v_{\varepsilon}}{1+v_{\varepsilon}m_{\varepsilon}},\left|\nabla m_{\varepsilon}\right|^2\right)\,dt=&4\int_0^T\left(\frac{v_{\varepsilon}m_{\varepsilon}}{1+v_{\varepsilon}m_{\varepsilon}},\left|\nabla m_{\varepsilon}^{\frac{1}{2}}\right|^2\right)\,dt\nonumber\\
 \leq&4\int_0^T\left(\frac{v_{\varepsilon}c_{\varepsilon}}{1+v_{\varepsilon}c_{\varepsilon}},\left|\nabla m_{\varepsilon}^{\frac{1}{2}}\right|^2\right)\,dt\nonumber\\
 \leq&\Cl{C2m}(T).\label{a4}
\end{align}
% {\cb Haben wir irgendwo gesagt, dass (und wieso) $0\le v_\varepsilon \le 1$?} 
Consequently, we also have that
\begin{align}
 \int_0^T\left(v_{\varepsilon},\left|\nabla\left(1+m_{\varepsilon}\right)^{\frac{1}{2}}\right|^2\right)\,dt=&\frac{1}{4}\int_0^T\left(\frac{v_{\varepsilon}}{1+m_{\varepsilon}},\left|\nabla m_{\varepsilon}\right|^2\right)\,dt\nonumber\\
 \leq&\frac{1}{4}\int_0^T\left(\frac{v_{\varepsilon}}{1+v_{\varepsilon}m_{\varepsilon}},\left|\nabla m_{\varepsilon}\right|^2\right)\,dt\nonumber\\
 \leq &\C(T).\label{a4_1}
\end{align}
\subsection*{Uniform integrability of $m_{\varepsilon}$}
It follows with \eqref{a1} that
\begin{align}
  &||m_{\varepsilon}||_{L^{\infty}(0,T;\LOne)}\leq\C(T).\label{c1}
\end{align}
Moreover, due to the de la Vall\'ee-Poussin theorem, we conclude with \eqref{a1} that
\begin{align}
 \{m_{\varepsilon}\}\text{ is uniformly integrable in }(0,T)\times\Omega.\label{mui}
\end{align}
\subsection*{Uniform integrability of  
 $\nabla\left(v_{\varepsilon}^{\frac{1}{2}}(m_{\varepsilon}+1)^{\frac{1}{2}}\right)$}
 Due to \eqref{mui}, it holds that
\begin{align}
 \left\{v_{\varepsilon}^{\frac{1}{2}}(m_{\varepsilon}+1)^{\frac{1}{2}}\right\}\text{ is uniformly integrable in }(0,T)\times\Omega.\label{d1ui}
\end{align}
We compute that
\begin{align}
 \nabla\left(v_{\varepsilon}^{\frac{1}{2}}(m_{\varepsilon}+1)^{\frac{1}{2}}\right)=v_{\varepsilon}^{\frac{1}{2}}\nabla(m_{\varepsilon}+1)^{\frac{1}{2}}+(m_{\varepsilon}+1)^{\frac{1}{2}}\nabla v_{\varepsilon}^{\frac{1}{2}}.\label{chmv1}
\end{align}
Combining \eqref{a2},  \eqref{a4_1}, \eqref{mui}, \eqref{chmv1} and using the de la Vall\'ee-Poussin theorem and {\it Lemma \ref{UiL2}}, we conclude that
\begin{align}
 \left\{\nabla\left(v_{\varepsilon}^{\frac{1}{2}}(m_{\varepsilon}+1)^{\frac{1}{2}}\right)\right\}\text{ is uniformly integrable in }(0,T)\times\Omega.\label{dui}
\end{align}
\subsection*{Uniform integrability of  the reaction term in \eqref{eq1e}}
It immediately follows with \eqref{cpb}, \eqref{c1}, \eqref{mui} that
\begin{align}
 \{f_{\varepsilon}\}\text{ is uniformly integrable in }(0,T)\times\Omega\label{fui}
\end{align}
and 
\begin{align}
 ||f_{\varepsilon}||_{L^{\infty}(0,T;\LOne)}\leq \C(T).\label{f1}
\end{align}
\subsection*{Uniform integrability of the diffusion flux in \eqref{eq1e}}
We first deal with the relaxation term. We have that 
\begin{align}
 \varepsilon_1^{\frac{1}{2}}|\nabla m_{\varepsilon}|&=2\varepsilon_1^{\frac{1}{2}}\left|\nabla m_{\varepsilon}^{\frac{1}{2}}\right|m_{\varepsilon}^{\frac{1}{2}},\label{este2}
\end{align}
Using the H\"older inequality, we obtain with \eqref{a6},  \eqref{c1} and \eqref{este2} that
\begin{align}
 \varepsilon_1||\nabla m_{\varepsilon}||_{L^1(0,T;\LOne)}\leq \varepsilon_1^{\frac{1}{2}}\C(T).\label{reld}
\end{align}
For the degenerate part of the diffusion flux, we have that 
\begin{align}
 \frac{v_{\varepsilon}c_{\varepsilon}|\nabla m_{\varepsilon}|}{1+v_{\varepsilon}c_{\varepsilon}}=&2\left(\frac{v_{\varepsilon}c_{\varepsilon}}{1+v_{\varepsilon}c_{\varepsilon}}\right)^{\frac{1}{2}}m_{\varepsilon}^{\frac{1}{2}}\left(\frac{v_{\varepsilon}c_{\varepsilon}}{1+v_{\varepsilon}c_{\varepsilon}}\left|\nabla m_{\varepsilon}^{\frac{1}{2}}\right|^2\right)^{\frac{1}{2}}.\label{dif1}
\end{align}
Combining \eqref{a4c} and \eqref{mui},
% and using the fact that $s\mapsto\frac{s}{1+s}$ is bounded
 we obtain with {\it Lemma \ref{UiL2}} that
\begin{align}
 \left\{\frac{v_{\varepsilon}c_{\varepsilon}\nabla m_{\varepsilon}}{1+v_{\varepsilon}c_{\varepsilon}}\right\}\text{ is uniformly integrable in }(0,T)\times\Omega,\label{ddui}
\end{align}
so that 
\begin{align}
 \left\|\frac{v_{\varepsilon}c_{\varepsilon}|\nabla m_{\varepsilon}|}{1+v_{\varepsilon}c_{\varepsilon}}\right\|_{L^1(0,T;\LOne)}\leq \C(T).\label{dflow1}
\end{align}
\subsection*{Uniform integrability of the taxis flux in \eqref{eq1e}}
Let us next consider the taxis part of the flux. We compute that
\begin{align}
 \frac{m_{\varepsilon}}{(1+v_{\varepsilon})^2}\nabla v_{\varepsilon}=&m_{\varepsilon}\nabla \frac{v_{\varepsilon}}{1+v_{\varepsilon}}\nonumber\\
 =&\nabla\left(\frac{v_{\varepsilon}m_{\varepsilon}}{1+v_{\varepsilon}}\right)-\frac{v_{\varepsilon}}{1+v_{\varepsilon}}\nabla m_{\varepsilon}.\label{taxdec}
\end{align}
For the second summand on the right-hand side of \eqref{taxdec}, we have that
\begin{align}
 \frac{v_{\varepsilon}}{1+v_{\varepsilon}}\nabla m_{\varepsilon}=2\frac{v_{\varepsilon}^{\frac{1}{2}}}{1+v_{\varepsilon}}(m_{\varepsilon}+1)^{\frac{1}{2}}v_{\varepsilon}^{\frac{1}{2}}\nabla(m_{\varepsilon}+1)^{\frac{1}{2}}\label{tax1}
\end{align}
We use \eqref{a4_1}, \eqref{mui} and {\it Lemma \ref{UiL2}} in order to conclude from \eqref{tax1} that
\begin{align}
 \left\{\frac{v_{\varepsilon}}{1+v_{\varepsilon}}\nabla m_{\varepsilon}\right\}\text{ is uniformly integrable in }(0,T)\times\Omega.\label{ta1ui}
\end{align}
As for the first summand on the right-hand side of \eqref{taxdec}, we seek for an estimate for its integral over $(0,t)$ (compare {\it Definition \ref{Defweak}}).
On both sides of equation \eqref{eq2e}, we divide by $1+v_{\varepsilon}$, apply the space gradient  and finally integrate over $(0,t)$. This yields that
\begin{align}
 \frac{1}{1+v_{\varepsilon}}\nabla v_{\varepsilon}(t)- \frac{1}{1+v_{\varepsilon_30}}\nabla v_{\varepsilon_30}=\mu_v\int_0^t\left(\frac{s(1-s)}{1+s}\right)'|_{s=v_{\varepsilon}}\nabla v_{\varepsilon}\,d\tau-\lambda\nabla\left(\int_0^t\frac{v_{\varepsilon}m_{\varepsilon}}{1+v_{\varepsilon}}\,d\tau\right).\label{vmod1}
\end{align}
Since $s\mapsto \frac{s(1-s)}{1+s}$ is continuously differentiable, we conclude from \eqref{vmod1} using \eqref{a2} that
\begin{align}
 \left\|\nabla\left(\int_0^t\frac{v_{\varepsilon}m_{\varepsilon}}{1+v_{\varepsilon}}\,d\tau\right)\right\|_{L^{\infty}(0,T;\LTwo)}\leq \C(T).
\end{align}
\subsection*{Estimates involving $\partial_t v_{\varepsilon}$}
We divide equation \eqref{eq2e} by $v_{\varepsilon}$:
\begin{align}
 \frac{1}{v_{\varepsilon}}\partial_t v_{\varepsilon}=\mu_v(1-v_{\varepsilon})-\lambda m_{\varepsilon}.\label{vmvt}
\end{align}
Together with  \eqref{c1}, \eqref{vmvt} yields that
\begin{align}
 \left\|\frac{1}{v_{\varepsilon}}\partial_t v_{\varepsilon}\right\|_{L^{\infty}(0,T;\LOne)}\leq \C(T),\label{lnvt}
\end{align}
so that, consequently, 
\begin{align}
 \left\|\partial_t v_{\varepsilon}^{\frac{1}{2}}\right\|_{L^{\infty}(0,T;\LOne)}\leq \C(T).\label{vt1}
\end{align}
\subsection*{Estimates for $\ln\left(1+v_{\varepsilon}^{\frac{1}{2}}m_{\varepsilon}\right)$}
%\texorpdfstring{$\ln\left(1+v_{\varepsilon}^{\frac{1}{2}}c_{\varepsilon}\right)$}{}}
Above we obtained uniform (in $\varepsilon$) estimates for both time and spacial derivatives of $v_{\varepsilon}$.
Owing to the fact that the original diffusion coefficient in \eqref{cm} is degenerate in $v$, it does not seem possible to obtain similar estimates for $m_{\varepsilon}$ or, at least, for $\varphi(m_{\varepsilon})$ for a     smooth, strictly increasing, and independent of $\varepsilon$ function  
$\varphi$. 
% Thus, one cannot apply a Lions-Aubin-type result directly to $\{m_{\varepsilon}\}$.
In order to overcome this difficulty and gain some information on $m_{\varepsilon}$ in the set  $\{v_{\varepsilon}>0\}$, we introduce for $\varepsilon\in(0,1)$ an auxiliary function which involves {\it both} $m_{\varepsilon}$ and $v_{\varepsilon}$:
\begin{align}
  u_{\varepsilon}:=\ln\left(1+v_{\varepsilon}^{\frac{1}{2}}m_{\varepsilon}\right). \label{u}%\nonumber%, \text{ where }  \varphi(s):=\ln(1+s).\label{u}
\end{align}
Since
\begin{align}
 0\leq\ln\left(1+v_{\varepsilon}^{\frac{1}{2}}m_{\varepsilon}\right)\leq m_{\varepsilon},\nonumber
\end{align}
we obtain with \eqref{mui} that
\begin{align}
 \{u_{\varepsilon}\}\text{ is uniformly integrable in }(0,T)\times\Omega.\nonumber
\end{align}
As it turns out, the family $\{u_{\varepsilon}\}$ is (strongly) precompact in $L^1(0,T;\LOne)$. To prove this, we need uniform estimates for the partial derivatives of $u_{\varepsilon}$ in some parabolic Sobolev spaces.

\noindent
We first study the spatial gradient of $u_{\varepsilon}$.
We compute that 
\begin{align}
 \nabla u_{\varepsilon}=\frac{m_{\varepsilon}}{1+v_{\varepsilon}^{\frac{1}{2}}m_{\varepsilon}}\nabla v_{\varepsilon}^{\frac{1}{2}}+\frac{v_{\varepsilon}^{\frac{1}{2}}}{1+v_{\varepsilon}^{\frac{1}{2}}m_{\varepsilon}}\nabla m_{\varepsilon}.\label{gru}
\end{align}
Using the trivial inequality 
\begin{align}
1\leq v_{\varepsilon}^{\frac{1}{2}}+(1-v_{\varepsilon})^{\frac{1}{2}},   \label{triv}                                           \end{align}
 we estimate the first summand on the right-hand side of \eqref{gru} in the following way:
\begin{align}
 \frac{m_{\varepsilon}\left|\nabla v_{\varepsilon}^{\frac{1}{2}}\right|}{1+v_{\varepsilon}^{\frac{1}{2}}m_{\varepsilon}}\leq&\frac{v_{\varepsilon}^{\frac{1}{2}}m_{\varepsilon}\left|\nabla v_{\varepsilon}^{\frac{1}{2}}\right|}{1+v_{\varepsilon}^{\frac{1}{2}}m_{\varepsilon}}+\frac{(1-v_{\varepsilon})^{\frac{1}{2}}m_{\varepsilon}\left|\nabla v_{\varepsilon}^{\frac{1}{2}}\right|}{1+v_{\varepsilon}^{\frac{1}{2}}m_{\varepsilon}}\nonumber\\
 \leq&\left|\nabla v_{\varepsilon}^{\frac{1}{2}}\right|+m_{\varepsilon}^{\frac{1}{2}}\left((1-v_{\varepsilon})m_{\varepsilon}\left|\nabla v_{\varepsilon}^{\frac{1}{2}}\right|^2\right)^{\frac{1}{2}}\nonumber\\
 \leq&\left|\nabla v_{\varepsilon}^{\frac{1}{2}}\right|+\frac{1}{2}m_{\varepsilon}+\frac{1}{2}(1-v_{\varepsilon})m_{\varepsilon}\left|\nabla v_{\varepsilon}^{\frac{1}{2}}\right|^2.
\label{estgr1}
\end{align}
Using estimates \eqref{a2}, \eqref{a5},  \eqref{c1}, we conclude from \eqref{estgr1} that
\begin{align}
 \left\|\frac{m_{\varepsilon}\left|\nabla v_{\varepsilon}^{\frac{1}{2}}\right|}{1+v_{\varepsilon}^{\frac{1}{2}}m_{\varepsilon}}\right\|_{L^1(0,T;\LOne)}%\nonumber\\
%  \leq& \C(T)\left\|\nabla v_{\varepsilon}^{\frac{1}{2}}\right\|_{L^2((0,T)\times\Omega)}+||m_{\varepsilon}||_{L^1(0,T;\LOne)}^{\frac12}\left(\int_0^T\left((1-v_{\varepsilon})m_{\varepsilon},\left|\nabla v_{\varepsilon}^{\frac{1}{2}}\right|^2\right)\,dt\right)^{\frac{1}{2}}\nonumber\\
 \leq&\C(T).\label{estgr1_}
\end{align}
For the second summand on the right-hand side of \eqref{gru}, we have that
\begin{align}
 &\frac{v_{\varepsilon}^{\frac{1}{2}}|\nabla m_{\varepsilon}|}{1+v_{\varepsilon}^{\frac{1}{2}}m_{\varepsilon}}\leq \frac{v_{\varepsilon}^{\frac{1}{2}}|\nabla m_{\varepsilon}|}{1+v_{\varepsilon}m_{\varepsilon}}\leq \frac{v_{\varepsilon}^{\frac{1}{2}}|\nabla m_{\varepsilon}|}{(1+v_{\varepsilon}m_{\varepsilon})^{\frac{1}{2}}}.\label{estgr2}
\end{align}
Due to \eqref{a4}, \eqref{estgr2} yields that
 \begin{align}
  \left\|\frac{v_{\varepsilon}^{\frac{1}{2}}|\nabla m_{\varepsilon}|}{1+v_{\varepsilon}^{\frac{1}{2}}m_{\varepsilon}}\right\|_{L^2((0,T)\times\Omega)}\leq\C(T).\label{estgr2_}
 \end{align}
Altogether, we obtain from \eqref{gru} with \eqref{estgr1_}, \eqref{estgr2_} that
\begin{align}
 ||\nabla u_{\varepsilon}||_{L^1(0,T;\LOne)}\leq\C(T).\label{estgru}
\end{align}
Next, we deal with the time derivative of $u_{\varepsilon}$. We compute that
\begin{align}
 \partial_t u_{\varepsilon}=\frac{1}{2}\frac{v_{\varepsilon}^{\frac{1}{2}}m_{\varepsilon}}{1+v_{\varepsilon}^{\frac{1}{2}}m_{\varepsilon}}\frac{1}{v_{\varepsilon}}\partial_t v_{\varepsilon}+\frac{v_{\varepsilon}^{\frac{1}{2}}}{1+v_{\varepsilon}^{\frac{1}{2}}m_{\varepsilon}}\partial_t m_{\varepsilon}.\label{tu}
\end{align}
We estimate the first summand on the right-hand side of \eqref{tu}  as follows:
\begin{align}
 \frac{1}{2}\left|\frac{v_{\varepsilon}^{\frac{1}{2}}m_{\varepsilon}}{1+v_{\varepsilon}^{\frac{1}{2}}m_{\varepsilon}}\frac{1}{v_{\varepsilon}}\partial_t v_{\varepsilon}\right|\leq \frac{1}{2}\left|\frac{1}{v_{\varepsilon}}\partial_t v_{\varepsilon}\right|.\label{est1s}
\end{align}
Combining \eqref{lnvt} and \eqref{est1s}, we obtain that
\begin{align}
 \frac{1}{2}\left\|\frac{v_{\varepsilon}^{\frac{1}{2}}m_{\varepsilon}}{1+v_{\varepsilon}^{\frac{1}{2}}m_{\varepsilon}}\frac{1}{v_{\varepsilon}}\partial_t v_{\varepsilon}\right\|_{L^{\infty}(0,T;\LOne)}\leq\C.\label{axx0}
\end{align}
In order to estimate the second summand on the right-hand side of \eqref{tu}, we multiply both sides of equation \eqref{eq1e} by $\frac{v_{\varepsilon}^{\frac{1}{2}}}{1+v_{\varepsilon}^{\frac{1}{2}}m_{\varepsilon}}$ and obtain (compare the notation in \eqref{qe}-\eqref{fe})  that
\begin{align}
 \frac{v_{\varepsilon}^{\frac{1}{2}}}{1+v_{\varepsilon}^{\frac{1}{2}}m_{\varepsilon}}\partial_t m_{\varepsilon}=\nabla\cdot \left(\frac{v_{\varepsilon}^{\frac{1}{2}}}{1+v_{\varepsilon}^{\frac{1}{2}}m_{\varepsilon}}q_{\varepsilon}\right)- q_{\varepsilon}\cdot\nabla \frac{v_{\varepsilon}^{\frac{1}{2}}}{1+v_{\varepsilon}^{\frac{1}{2}}m_{\varepsilon}} +\frac{v_{\varepsilon}^{\frac{1}{2}}}{1+v_{\varepsilon}^{\frac{1}{2}}m_{\varepsilon}}f_{\varepsilon}.\label{ax1}
\end{align}
Since 
\begin{align}
 \frac{v_{\varepsilon}^{\frac{1}{2}}}{1+v_{\varepsilon}^{\frac{1}{2}}m_{\varepsilon}}\leq1,\nonumber
\end{align}
% estimates \eqref{f1} and \eqref{dif2} yield, respectively, that
it holds that 
\begin{align}
 \frac{v_{\varepsilon}^{\frac{1}{2}}}{1+v_{\varepsilon}^{\frac{1}{2}}m_{\varepsilon}}|f_{\varepsilon}|\leq|f_{\varepsilon}|.\nonumber
 \end{align}
 Hence, we conclude with \eqref{f1} that
 \begin{align}
  \left\|\frac{v_{\varepsilon}^{\frac{1}{2}}}{1+v_{\varepsilon}^{\frac{1}{2}}m_{\varepsilon}}f_{\varepsilon}\right\|_{L^{\infty}(0,T;\LOne)}\leq \C(T).\label{axf}
 \end{align}
 For the term inside the divergence operator in \eqref{ax1}, we have that
\begin{align}
 \frac{v_{\varepsilon}^{\frac{1}{2}}}{1+v_{\varepsilon}^{\frac{1}{2}}m_{\varepsilon}}|q_{\varepsilon}|\leq&\frac{v_{\varepsilon}^{\frac{1}{2}}}{1+v_{\varepsilon}^{\frac{1}{2}}m_{\varepsilon}}\left(\varepsilon_1|\nabla m_{\varepsilon}|+\kappa _m\frac{v_{\varepsilon}c_{\varepsilon}|\nabla m_{\varepsilon}|}{1+v_{\varepsilon}c_{\varepsilon}}+\kappa _v\frac{2 v_{\varepsilon}^{\frac{1}{2}}m_{\varepsilon}\left|\nabla v_{\varepsilon}^{\frac{1}{2}}\right|}{(1+v_{\varepsilon})^2}\right)\nonumber\\
 \leq&\varepsilon_1|\nabla m_{\varepsilon}|+\kappa _m\frac{v_{\varepsilon}c_{\varepsilon}|\nabla m_{\varepsilon}|}{1+v_{\varepsilon}c_{\varepsilon}}+2\kappa _v\left|\nabla v_{\varepsilon}^{\frac{1}{2}}\right|.\label{axx1}
\end{align}
Using \eqref{a2}, \eqref{reld}, \eqref{dflow1}, we obtain from \eqref{axx1} that
\begin{align}
 \left\|\frac{v_{\varepsilon}^{\frac{1}{2}}}{1+v_{\varepsilon}^{\frac{1}{2}}m_{\varepsilon}}q_{\varepsilon}\right\|_{L^1(0,T;\LOne)}\leq \C(T).\label{ax0}
\end{align}
It remains to estimate the second term on the right-hand side of \eqref{ax1}.
We compute that
\begin{align}
 \nabla\frac{v_{\varepsilon}^{\frac{1}{2}}}{1+v_{\varepsilon}^{\frac{1}{2}}m_{\varepsilon}}=-\frac{v_{\varepsilon}}{\left(1+v_{\varepsilon}^{\frac{1}{2}}m_{\varepsilon}\right)^2}\nabla m_{\varepsilon}+\frac{1}{\left(1+v_{\varepsilon}^{\frac{1}{2}}m_{\varepsilon}\right)^2}\nabla v_{\varepsilon}^{\frac{1}{2}},\nonumber
\end{align}
so that
\begin{align}
 \left|q_{\varepsilon}\cdot\nabla \frac{v_{\varepsilon}^{\frac{1}{2}}}{1+v_{\varepsilon}^{\frac{1}{2}}m_{\varepsilon}}\right|% & \nonumber\\
 &\leq|q_{\varepsilon}|\left|\nabla \frac{v_{\varepsilon}^{\frac{1}{2}}}{1+v_{\varepsilon}^{\frac{1}{2}}m_{\varepsilon}}\right|\nonumber\\
 \leq&\left(\varepsilon_1|\nabla m_{\varepsilon}|+\frac{\kappa_m v_{\varepsilon}c_{\varepsilon}|\nabla m_{\varepsilon}|}{1+v_{\varepsilon}c_{\varepsilon}}+\frac{2 v_{\varepsilon}^{\frac{1}{2}}m_{\varepsilon}\left|\nabla v_{\varepsilon}^{\frac{1}{2}}\right|}{(1+v_{\varepsilon})^2}\right)\left(\frac{v_{\varepsilon}|\nabla m_{\varepsilon}|}{\left(1+v_{\varepsilon}^{\frac{1}{2}}m_{\varepsilon}\right)^2}+\frac{\left|\nabla v_{\varepsilon}^{\frac{1}{2}}\right|}{\left(1+v_{\varepsilon}^{\frac{1}{2}}m_{\varepsilon}\right)^2}\right)
 \nonumber\\
 \leq&\C\left(\varepsilon_1|\nabla m_{\varepsilon}|+\frac{ v_{\varepsilon}c_{\varepsilon}|\nabla m_{\varepsilon}|}{1+v_{\varepsilon}c_{\varepsilon}}+v_{\varepsilon}^{\frac{1}{2}}m_{\varepsilon}\left|\nabla v_{\varepsilon}^{\frac{1}{2}}\right|\right)\left(\frac{v_{\varepsilon}|\nabla m_{\varepsilon}|}{\left(1+v_{\varepsilon}^{\frac{1}{2}}m_{\varepsilon}\right)^2}+\frac{\left|\nabla v_{\varepsilon}^{\frac{1}{2}}\right|}{\left(1+v_{\varepsilon}^{\frac{1}{2}}m_{\varepsilon}\right)^2}\right).\label{ax2}
\end{align}
Using \eqref{triv} { and \eqref{cpb}} where necessary, we get the following estimates:
\begin{align}
& |\nabla m_{\varepsilon}|\frac{v_{\varepsilon}|\nabla m_{\varepsilon}|}{\left(1+v_{\varepsilon}^{\frac{1}{2}}m_{\varepsilon}\right)^2}\leq \frac{v_{\varepsilon}|\nabla m_{\varepsilon}|^2}{1+v_{\varepsilon}m_{\varepsilon}},\label{ax3}
\end{align}
\begin{align}
& \frac{ v_{\varepsilon}c_{\varepsilon}|\nabla m_{\varepsilon}|}{1+v_{\varepsilon}c_{\varepsilon}}\frac{v_{\varepsilon}|\nabla m_{\varepsilon}|}{\left(1+v_{\varepsilon}^{\frac{1}{2}}m_{\varepsilon}\right)^2}\leq \frac{v_{\varepsilon}|\nabla m_{\varepsilon}|^2}{1+v_{\varepsilon}m_{\varepsilon}}
\end{align}
\begin{align}
 v_{\varepsilon}^{\frac{1}{2}}m_{\varepsilon}\left|\nabla v_{\varepsilon}^{\frac{1}{2}}\right|\frac{v_{\varepsilon}|\nabla m_{\varepsilon}|}{\left(1+v_{\varepsilon}^{\frac{1}{2}}m_{\varepsilon}\right)^2}\leq &\left(\frac{v_{\varepsilon}|\nabla m_{\varepsilon}|^2}{1+v_{\varepsilon}m_{\varepsilon}}\right)^{\frac{1}{2}}\left|\nabla v_{\varepsilon}^{\frac{1}{2}}\right|\nonumber\\
\leq &\frac{1}{2}\frac{v_{\varepsilon}|\nabla m_{\varepsilon}|^2}{1+v_{\varepsilon}m_{\varepsilon}}+\frac{1}{2}\left|\nabla v_{\varepsilon}^{\frac{1}{2}}\right|^2,
\end{align}
\begin{align}
 |\nabla m_{\varepsilon}|\frac{\left|\nabla v_{\varepsilon}^{\frac{1}{2}}\right|}{\left(1+v_{\varepsilon}^{\frac{1}{2}}m_{\varepsilon}\right)^2}\leq&\frac{v_{\varepsilon}^{\frac{1}{2}}|\nabla m_{\varepsilon}|\left|\nabla v_{\varepsilon}^{\frac{1}{2}}\right|}{\left(1+v_{\varepsilon}m_{\varepsilon}\right)^2}+\frac{(1-v_{\varepsilon})^{\frac{1}{2}}|\nabla m_{\varepsilon}|\left|\nabla v_{\varepsilon}^{\frac{1}{2}}\right|}{\left(1+v_{\varepsilon}m_{\varepsilon}\right)^2}\nonumber\\
 \leq &\left(\frac{v_{\varepsilon}|\nabla m_{\varepsilon}|^2}{1+v_{\varepsilon}m_{\varepsilon}}\right)^{\frac{1}{2}}\left|\nabla v_{\varepsilon}^{\frac{1}{2}}\right|+2\left|\nabla m_{\varepsilon}^{\frac{1}{2}}\right|\left((1-v_{\varepsilon})m_{\varepsilon}\left|\nabla v_{\varepsilon}^{\frac{1}{2}}\right|^2\right)^{\frac{1}{2}}\nonumber\\
 \leq&\frac{1}{2}\frac{v_{\varepsilon}|\nabla m_{\varepsilon}|^2}{1+v_{\varepsilon}m_{\varepsilon}}+\frac{1}{2}\left|\nabla v_{\varepsilon}^{\frac{1}{2}}\right|^2+\left|\nabla m_{\varepsilon}^{\frac{1}{2}}\right|^2+(1-v_{\varepsilon})m_{\varepsilon}\left|\nabla v_{\varepsilon}^{\frac{1}{2}}\right|^2,
\end{align}
{
\begin{align}
 \frac{v_{\varepsilon}c_{\varepsilon}|\nabla m_{\varepsilon}|}{1+v_{\varepsilon}c_{\varepsilon}}\frac{\left|\nabla v_{\varepsilon}^{\frac{1}{2}}\right|}{\left(1+v_{\varepsilon}^{\frac{1}{2}}m_{\varepsilon}\right)^2}\leq &\frac{2v_{\varepsilon}^{\frac{1}{2}}m_{\varepsilon}^{\frac{1}{2}}\left(m_{\varepsilon}+p_{\varepsilon}\right)^{\frac{1}{2}}}{\left(1+v_{\varepsilon}^{\frac{1}{2}}m_{\varepsilon}\right)\left(1+v_{\varepsilon}c_{\varepsilon}\right)^{\frac{1}{2}}}\left(\frac{v_{\varepsilon}c_{\varepsilon}\left|\nabla m_{\varepsilon}^{\frac{1}{2}}\right|^2}{1+v_{\varepsilon}c_{\varepsilon}}\right)^{\frac{1}{2}}\left|\nabla v_{\varepsilon}^{\frac{1}{2}}\right|\nonumber\\
 \leq&\frac{2v_{\varepsilon}^{\frac{1}{2}}m_{\varepsilon}+2C_p^{\frac{1}{2}}(v_{\varepsilon}m_{\varepsilon})^{\frac{1}{2}}}{\left(1+v_{\varepsilon}^{\frac{1}{2}}m_{\varepsilon}\right)\left(1+v_{\varepsilon}m_{\varepsilon}\right)^{\frac{1}{2}}}\left(\frac{v_{\varepsilon}c_{\varepsilon}\left|\nabla m_{\varepsilon}^{\frac{1}{2}}\right|^2}{1+v_{\varepsilon}c_{\varepsilon}}\right)^{\frac{1}{2}}\left|\nabla v_{\varepsilon}^{\frac{1}{2}}\right|\nonumber\\
 %{\cb \text{Was ist mit } 2\meps^{1/2} \text{ passiert?}}
\leq &\Cl{C100}\frac{v_{\varepsilon}c_{\varepsilon}\left|\nabla m_{\varepsilon}^{\frac{1}{2}}\right|^2}{1+v_{\varepsilon}c_{\varepsilon}}+\Cr{C100}\left|\nabla v_{\varepsilon}^{\frac{1}{2}}\right|^2,
\end{align}
\begin{align}
 v_{\varepsilon}^{\frac{1}{2}}m_{\varepsilon}\left|\nabla v_{\varepsilon}^{\frac{1}{2}}\right|\frac{\left|\nabla v_{\varepsilon}^{\frac{1}{2}}\right|}{\left(1+v_{\varepsilon}^{\frac{1}{2}}m_{\varepsilon}\right)^2}\leq &\left|\nabla v_{\varepsilon}^{\frac{1}{2}}\right|^2.\label{ax4}
\end{align}
}
% {\cb 
% \begin{align}
% \frac{\veps \ceps \ |\nabla \meps |}{1+\veps\ceps}\ \frac{\veps |\nabla \meps |}{(1+\veps ^{1/2}\meps )^2} &\le 2\Big (\frac{\veps \ceps \ |\nabla \meps ^{1/2}|^2}{1+\veps\ceps}\Big )^{1/2}\Big (\frac{\veps \ceps }{1+\veps\ceps}\Big )^{1/2}\meps ^{1/2}\frac{\veps |\nabla \meps |}{1+\veps ^{1/2}\meps }\nonumber \\
% &\le 2|\nabla \meps ^{1/2}|^2+2\frac{\veps \ceps \ |\nabla \meps ^{1/2}|^2}{1+\veps\ceps}.\label{ax4b}
% \end{align}
% }
\noindent
Combining \eqref{ax2}-\eqref{ax4} with \eqref{a2}-\eqref{a5}, \eqref{a6}, \eqref{a4}, \eqref{c1}, we obtain that
\begin{align}
 \left\|q_{\varepsilon}\cdot\nabla \frac{v_{\varepsilon}^{\frac{1}{2}}}{1+v_{\varepsilon}^{\frac{1}{2}}m_{\varepsilon}}\right\|_{L^1(0,T;\LOne)}\leq\C(T).\label{ax5}
\end{align}
Therefore,  \eqref{ax1} together with \eqref{axf}, \eqref{ax0} and \eqref{ax5} yield that
\begin{align}
 \left\|\frac{v_{\varepsilon}^{\frac{1}{2}}\partial_t m_{\varepsilon}}{1+v_{\varepsilon}^{\frac{1}{2}}m_{\varepsilon}}\right\|_{L^1(0,T;W^{-1,1}(\Omega))}\leq\C(T).\label{ax6}
\end{align}
Finally, with the help of estimates \eqref{axx0} and \eqref{ax6}, we obtain from \eqref{tu} that
\begin{align}
||\partial_t u_{\varepsilon}||_{L^1(0,T;W^{-1,1}(\Omega))}\leq \C(T).\label{uet}
\end{align}
\subsection*{Estimates for $m_{\varepsilon}$ in $(0,T)\times\inner\{v_{\varepsilon_30}=0\}$}
While studying the function $m_{\varepsilon}$, the auxiliary function $u_{\varepsilon}$ introduced in \eqref{u} is of use only in the set $\{v_{\varepsilon}>0\}$. It clearly reveals no further information about the behaviour of $m_{\varepsilon}$ over the level sets $\{v_{\varepsilon}(t)=0\}$, $t\in(0,T)$, each of whom almost coincide 
with $\{v_{\varepsilon_30}=0\}$. 
% Before we begin with the estimates, we make the following important observation:  equation \eqref{eq2e} is an ODE for $v_{\varepsilon}$, it preserves the null set of $v_{\varepsilon}$ with time. 
The latter is to mean that $\{v_{\varepsilon}(t)=0\}$ differs from  $\{v_{\varepsilon_30}=0\}$ by a null set and is thus preserved in time. In order to see this, let us divide both sides of the ODE \eqref{eq2e} by $v_{\varepsilon}$ and integrate over $(0,t)$ for arbitrary $t\in(0,T)$. We obtain that
\begin{align}
 \ln(v_{\varepsilon}(t))-\ln(v_{\varepsilon_30})=&\int_{0}^t\mu_v(1-v_{\varepsilon})\,dt-\lambda \int_{0}^tm_{\varepsilon}\,dt.\label{lnv}
\end{align}
Since $0\leq v_{\varepsilon}\leq1$ and $m_{\varepsilon}\in L^1(0,T;\LOne)$, the right-hand side of \eqref{lnv} is finite a.e. in $\Omega$. Hence, the same holds for the left-hand side of \eqref{lnv}. But this means that for all $t\in\R^+$ it necessarily holds that 
\begin{align}
   &v_{\varepsilon}(t)>0 \text{ a.e. in }\{v_{\varepsilon_30}>0\},\nonumber\\%\label{vposn}\\
   &v_{\varepsilon}(t)=0 \text{ a.e. in }\{v_{\varepsilon_30}=0\}.\label{vnposn}
\end{align}
Similarly, we obtain from the original equation \eqref{v} that
\begin{align}
   &v(t)>0 \text{ a.e. in }\{v_0>0\},\nonumber\\%\label{vpos}\\
   &v(t)=0 \text{ a.e. in }\{v_0=0\}.\label{vnpos}
\end{align}
Observe that, at least in  $(0,T)\times\inner\{v_{\varepsilon_30}=0\}$, $m_{\varepsilon}$ solves the linear initial value problem
\begin{subequations}\label{linm}
\begin{alignat}{3}
 & \partial_t m_{\varepsilon}=\varepsilon_1\Delta m_{\varepsilon}-\alpha m_{\varepsilon}&&\text{ in }\R^+\times \inner\{v_{\varepsilon_30}=0\},\label{mv0}\\
 & m_{\varepsilon}(0)=m_{\varepsilon_20}&&\text{ in }\inner\{v_{\varepsilon_30}=0\}.\label{eqmv0}
\end{alignat}
\end{subequations}
Combining \eqref{c1} and \eqref{reld}, we conclude from \eqref{mv0} that
\begin{align}
||\partial_t m_{\varepsilon}||_{L^1(0,T;W^{-1,1}(\inner\{v_{\varepsilon_30}=0\}))}\leq \C(T).\label{mv03t}
\end{align}
Since $m_{\varepsilon_20}$ is smooth, $m_{\varepsilon}$ is a classical solution to \eqref{mv0}. Differentiating \eqref{mv0} with respect to $x_i$, $i\in \{1,\dots ,N\}$, we obtain that
\begin{align}
& \partial_t \partial_{x_i} m_{\varepsilon}=\varepsilon_1\Delta\partial_{x_i}  m_{\varepsilon}-\alpha \partial_{x_i} m_{\varepsilon}.\label{nmv0}
\end{align}
Let now $\varphi$ be some smooth cut-off function with $\supp\varphi\subset \inner\{v_{\varepsilon_30}=0\}$ and let $a\in(1,2)$, the latter to be specified below. Multiplying \eqref{nmv0} by $a\varphi^2|\partial_{x_i} m_{\varepsilon}|^{a-2}\partial_{x_i} m_{\varepsilon}$ and integrating by parts over $\Omega$, we obtain with the H\"older and Young inequalities that
\begin{align}
 \frac{d}{dt}\left\|\varphi |\partial_{x_i} m_{\varepsilon}|^{\frac{a}{2}}\right\|^2
 =&-\frac{4(a-1)}{a}\varepsilon_1\left\|\varphi\nabla|\partial_{x_i} m_{\varepsilon}|^{\frac{a}{2}}\right\|^2-4\varepsilon_1\left(\varphi\nabla|\partial_{x_i} m_{\varepsilon}|^{\frac{a}{2}},|\partial_{x_i} m_{\varepsilon}|^{\frac{a}{2}}\nabla\varphi\right) -a\alpha\left\|\varphi |\partial_{x_i} m_{\varepsilon}|^{\frac{a}{2}}\right\|^2\nonumber\\
 \leq&\C(a)||\nabla\varphi||_{\infty}^2\varepsilon_1\left\|\partial_{x_i} m_{\varepsilon}\right\|^{a}_{a}\nonumber\\
 =&\C(a)||\nabla\varphi||_{\infty}^2\varepsilon_1\left\|m_{\varepsilon}^{\frac{1}{2}}\left|\partial_{x_i} m_{\varepsilon}^{\frac{1}{2}}\right|\right\|^{a}_{a}\nonumber\\
 \leq&\C(a)||\nabla\varphi||_{\infty}^2\varepsilon_1\left\|m_{\varepsilon}^{\frac{1}{2}}\right\|^{a}_{\frac{1}{\frac{1}{a}-\frac{1}{2}}}\left\|\partial_{x_i} m_{\varepsilon}^{\frac{1}{2}}\right\|^{a}.\label{estmv01}
\end{align}
Owing to a Sobolev interpolation inequality it holds that
\begin{align}
 \left\|m_{\varepsilon}^{\frac{1}{2}}\right\|_{\frac{1}{\frac{1}{a}-\frac{1}{2}}}\leq \C(a)\left(\left\|\nabla m_{\varepsilon}^{\frac{1}{2}}\right\|+\left\|m_{\varepsilon}^{\frac{1}{2}}\right\|\right)^{N\left(1-\frac{1}{a}\right)}\left\|m_{\varepsilon}^{\frac{1}{2}}\right\|^{1-N\left(1-\frac{1}{a}\right)}\text{ for }a\in\left[1,\frac{N}{N-1}\right].\label{estax}
\end{align}
Integrating \eqref{estmv01} over $(0,t)$ and using \eqref{a6}, \eqref{c1}, \eqref{estax} and the H\"older inequality, we thus obtain that
\begin{align}
 \left\|\varphi |\nabla m_{\varepsilon}(t)|^{\frac{a}{2}}\right\|^2\leq&\left\|\varphi |\nabla m_{\varepsilon_20}|^{\frac{a}{2}}\right\|^2+\C(a,T)||\nabla\varphi||_{\infty}^2\varepsilon_1^{1-\frac{N}{2}(a-1)-\frac{a}{2}}\text{ for }a\in\left(1,\frac{N}{N-1}\right].\label{mv0_}
\end{align} 
The first term on right-hand side of \eqref{mv0_} doesn't depend upon $\varepsilon_1$, while the second one converges to zero for $a\in\left(1,\frac{N+2}{N+1}\right)\subset\left(1,\frac{N}{N-1}\right]$. Therefore, we obtain from \eqref{mv0_} in the limit as $\varepsilon_1\rightarrow0$ that
\begin{align}
 \underset{\varepsilon_1\rightarrow0}{\lim\sup}\left\|\varphi |\nabla m_{\varepsilon}(t)|^{\frac{a}{2}}\right\|^2\leq&\left\|\varphi |\nabla m_{\varepsilon_20}|^{\frac{a}{2}}\right\|^2\nonumber\\
 \leq&||\nabla m_{\varepsilon_20}||^a_{\infty}\text{ for }a\in\left(1,\frac{N+2}{N+1}\right).\label{mv01}
\end{align}
Since $\varphi$ was an arbitrary cut-off function with $\supp\varphi\subset \inner\{v_{\varepsilon_30}=0\}$, \eqref{mv01} yields that 
\begin{align}
 \underset{\varepsilon_1\rightarrow0}{\lim\sup}||\nabla m_{\varepsilon}(t)||_{L^a(\inner\{v_{\varepsilon_30}=0\})}
 \leq&||\nabla m_{\varepsilon_20}||_{\infty}\text{ for }a\in\left(1,\frac{N+2}{N+1}\right).\label{mv02_}
\end{align}
Together with \eqref{c1}, \eqref{mv02_} yields that
\begin{align}
 \underset{\varepsilon_1\rightarrow0}{\lim\sup}||m_{\varepsilon}(t)||_{W^{1,1}(\inner\{v_{\varepsilon_30}=0\})}
 \leq&\C(\varepsilon_2).\label{mv03}
\end{align}
\section{Global existence for the original problem}\label{existence}
In this section we aim to pass to the limit in \eqref{hapto2e} in order to obtain a solution of the original problem. 
\begin{Remark}[Vector notation]
 Let $ \{\varepsilon_{i,n_{i}}\}\subset (0,1)$, $i= 1,2,3$, be three sequences.
 %such that for each $i= 1,2,3$ it holds that $\varepsilon_{i,n_{i}}
% \underset{n_{i}\rightarrow\infty}{\rightarrow}0$. 
In this section, we make use of the following vector notation:
\begin{align}
 &n_{i:3}:=\left(n_{i},\dots,n_3\right), \
\varepsilon_{n_{i:3}}:=\left(\varepsilon_{i,n_{i}},\dots,\varepsilon_{3,n_3}\right),\ i=1,2.\nonumber
\end{align}
% Let us illustrate the way we are going to apply it. Let $a_{\varepsilon_{n_{1:3}}}$ be a family parameterized by $\varepsilon$, i.e., by the triple $\left(\varepsilon_{1},\varepsilon_{2},\varepsilon_3\right)$. By writing
% \begin{align}
%  a_{\varepsilon_{n_{1:3}}}\underset{n_{1}\rightarrow\infty}{\rightarrow}a_{n_{2:3}}\underset{n_{1:3}\rightarrow\infty}{\rightarrow} a\text{ (in some topology)},\label{limit}
% \end{align}
% we mean thus that the sequence $\{a_{\varepsilon_{n_{1:3}}}\}$ converges to some  $a_{n_{2:3}}$ as $n_{1}\rightarrow\infty$ for each $n_{2:3}$, while $\{a_{n_{2:3}}\}$ converges 
% to some $a$ (in some topology) as $n_{2:3}\rightarrow\infty$, i.e., as $n_{2},n_{3}\rightarrow\infty$. 
% As for the family of limits $\{a_{n_{2:3}}\}$ and the limit $a$, it is assumed that they  either have been previously introduced, or that they exist and are  being thus 
% introduced by expression \eqref{limit}. \\
% \noindent
% Thereby, we can write two subsequent limit procedures in a compact form.
% % Thus, $a$ is a repeated limit of $\{a_{\varepsilon_{n_{1:3}}}\}$:
% % \begin{align}
% %  a_{\varepsilon_{n_{1:3}}}\underset{n_{1}\rightarrow\infty}{\rightarrow}\underset{n_{1:3}\rightarrow\infty}{\rightarrow} a\text{ in [some sense]}.\nonumber
% % \end{align}
% 
\end{Remark}
\noindent
% From know on, we assume that the families $\{m_{\varepsilon_20}\}$, $\{p_{\varepsilon_20}\}$ and $\{v_{\varepsilon_30}\}$ are parameterized by $\varepsilon_2$ and $\varepsilon_3$, respectively, and   satisfy \eqref{e34bnd}-\eqref{star}. 
Owing to the estimates obtained in the preceding section, we are now in a position where we can establish a list convergences (see below) holding jointly for some sequences $$ \varepsilon_{i,n_{i}}\underset{n_{i}\rightarrow\infty}{\rightarrow}0,\ i= 1,2,3.$$ 
\subsection*{Convergence for the initial data}
Due to  \eqref{moe}-\eqref{voe} it holds that
\begin{align}
 &m_{\varepsilon_{2,n_2}0}\underset{n_{2}\rightarrow\infty}{\rightarrow}m_0\text{ in }L^1(\Omega)\text{ and a.e. in }\Omega,\label{mk0}\\
 &p_{\varepsilon_{2,n_2}0}\underset{n_{2}\rightarrow\infty}{\rightarrow}p_0\text{ in }L^{\infty}(\Omega)\text{ and a.e. in }\Omega,\label{pk0}\\
 &v_{\varepsilon_{3,n_3}0}^{\frac{1}{2}}\underset{n_3\rightarrow\infty}{\rightarrow}v^{\frac{1}{2}}_0\text{ in }\LTwo\text{ and a.e. in }\Omega;\label{vk0}
\end{align}
\subsection*{Convergence for $\left\{v_{\varepsilon_{n_{1:3}}}\right\}$}It holds that: 
due to \eqref{a2}, \eqref{vt1} and a version of the Lions-Aubin Lemma \cite[Corollary 4]{Simon} 
\begin{align}
  v_{\varepsilon_{n_{1:3}}}^{\frac{1}{2}}\underset{n_{1:3}\rightarrow\infty}{\rightarrow} v^{\frac{1}{2}}\text{ in }L^2(0,T;\LTwo);\label{vcomp}
\end{align}
due to \eqref{vcomp}
\begin{align}
 v_{\varepsilon_{n_{1:3}}}^{\frac{1}{2}}\underset{n_{1:3}\rightarrow\infty}{\rightarrow} v^{\frac{1}{2}}\text{ a.e. in  }(0,T)\times\Omega;\label{vae}
\end{align}
due to \eqref{vae}
\begin{align}
 v_{\varepsilon_{n_{1:3}}}\underset{n_{1:3}\rightarrow\infty}{\rightarrow} v\text{ a.e. in  }(0,T)\times\Omega;\label{vae_}
\end{align}
due to  \eqref{vae_} and the  dominated convergence theorem
\begin{align}
 v_{\varepsilon_{n_{1:3}}}^a\underset{n_{1:3}\rightarrow\infty}{\rightarrow} v^a\text{ in }L^{p}((0,T)\times\Omega)\text{ and a.e. in }(0,T)\times\Omega\text{ for all }a>0,\ p\geq1;\label{vconvp}
\end{align}
due to \eqref{a2}, \eqref{vcomp} and the Banach-Alaoglu theorem
\begin{align}
 \nabla v_{\varepsilon_{n_{1:3}}}^{\frac{1}{2}}\underset{n_{1:3}\rightarrow\infty}{\rightharpoonup}\nabla v^{\frac{1}{2}}\text{ in }L^2(0,T;\LTwo).\label{nv2}
\end{align}
\subsection*{Convergence for $\left\{m_{\varepsilon_{n_{1:3}}}\right\}$ in $(0,T)\times\{v_0>0\}$}
% due to \eqref{ui} and the Dunford-Pettis theorem
% \begin{align}
% &m_{\varepsilon_{n_{1:3}}}\underset{n_{1:3}\rightarrow\infty}{\rightharpoonup}m\text{ in }L^1(0,T;\LOne);\label{cmconv}
% \end{align}
It holds that: due to \eqref{estgru}, \eqref{uet} and a version of the Lions-Aubin Lemma \cite[Corollary 4]{Simon}
\begin{align}
 \ln\left(1+v_{\varepsilon_{n_{1:3}}}^{\frac{1}{2}}m_{\varepsilon_{n_{1:3}}}\right)\underset{n_{1:3}\rightarrow\infty}{\rightarrow} u\text{ in }L^1(0,T;\LOne);\label{ucomp2}
\end{align}
due to \eqref{ucomp2}
\begin{align}
 \ln\left(1+v_{\varepsilon_{n_{1:3}}}^{\frac{1}{2}}m_{\varepsilon_{n_{1:3}}}\right)\underset{n_{1:3}\rightarrow\infty}{\rightarrow} u\text{ a.e. in  }(0,T)\times\Omega;\label{uae}
\end{align}
due to \eqref{uae}
\begin{align}
 v_{\varepsilon_{n_{1:3}}}^{\frac{1}{2}}m_{\varepsilon_{n_{1:3}}}\underset{n_{1:3}\rightarrow\infty}{\rightarrow} e^u-1=:w\text{ a.e. in  }(0,T)\times\Omega;\label{wae}
\end{align}
due to \eqref{vae_}, \eqref{wae}% and the Lions lemma \cite[Lemma 1.3]{Lions} (continues to hold in $L^1$...)
\begin{align}
 m_{\varepsilon_{n_{1:3}}}\underset{n_{1:3}\rightarrow\infty}{\rightarrow}\frac{w}{v^{\frac{1}{2}}}=:m\text{ a.e. in  }(0,T)\times\{v_0>0\};\label{mae}
\end{align}
due to \eqref{mui}, \eqref{mae} and the Vitali convergence theorem
\begin{align}
 m_{\varepsilon_{n_{1:3}}}\underset{n_{1:3}\rightarrow\infty}{\rightarrow}m\text{ in  }L^1((0,T)\times\{v_0>0\}).\label{cmstr}
\end{align}
\subsection*{Convergence for $\left\{m_{\varepsilon_{n_{1:3}}}\right\}$ in $(0,T)\times\{v_0=0\}$}
It holds due to \eqref{mv03t}, \eqref{mv03} and a version of the Lions-Aubin Lemma \cite[Corollary 4]{Simon} that
\begin{align}
 m_{\varepsilon_{n_{1:3}}}\underset{n_1\rightarrow\infty}{\rightarrow}m_{n_{2:3}}\text{ in  }L^1((0,T)\times\inner\{v_{\varepsilon_{3,n_3}0}=0\}),\label{m123conv}
\end{align}
and so we may pass to the distributional limit in \eqref{linm}:
\begin{subequations}\label{linm2}
\begin{alignat}{3}
 & \partial_t m_{n_{2:3}}=-\alpha m_{n_{2:3}}&&\text{ in }\R^+\times \inner\{v_{\varepsilon_{3,n_3}0}=0\},\label{mv02}\\
 & m_{n_{2:3}}(0)=m_{\varepsilon_{2,n_2}0}&&\text{ in }\inner\{v_{\varepsilon_{3,n_3}0}=0\}.\label{eqmv02}
\end{alignat}
\end{subequations}
Due to \eqref{mk0} and the continuous dependence of solutions of an ODE with smooth coefficients upon the initial data, it follows with \eqref{eqmv02} that
\begin{align}
 m_{n_{2:3}}\underset{n_2\rightarrow\infty}{\rightarrow}m_{n_{3}}\text{ in  }L^1((0,T)\times\inner\{v_{\varepsilon_{3,n_3}0}=0\})\label{m23conv}
\end{align}
and 
\begin{align}
 m_{n_{3}}=m\text{ a.e. in }(0,T)\times\inner\{v_{\varepsilon_{3,n_3}0}=0\}\cap\{v_0=0\},
\end{align}
where $m$ solves
\begin{subequations}\label{linm3}
\begin{alignat}{3}
 & \partial_t m=-\alpha m&&\text{ in }\R^+\times \{v_0=0\},\\
 & m(0)=m_0&&\text{ in }\{v_0=0\}.
\end{alignat}
\end{subequations}
Combining \eqref{m123conv}, \eqref{m23conv}-\eqref{linm3}, we conclude that
\begin{align}
 m_{\varepsilon_{n_{1:3}}}\underset{n_1\rightarrow\infty}{\rightarrow}\underset{n_2\rightarrow\infty}{\rightarrow}m\text{ in  }L^1((0,T)\times\inner\{v_{\varepsilon_{3,n_3}0}=0\}\cap\{v_0=0\}),
\end{align}
hence also
\begin{align}
 m_{\varepsilon_{n_{1:3}}}\underset{n_1\rightarrow\infty}{\rightarrow}\underset{n_2\rightarrow\infty}{\rightarrow}m\text{ on  }(0,T)\times\inner\{v_{\varepsilon_{3,n_3}0}=0\}\cap\{v_0=0\}\text{ in the measure}.\label{m1-3conv}
\end{align}
Together with property \eqref{star}, \eqref{m1-3conv} yields that 
\begin{align}
 \underset{n_3\rightarrow\infty}{\lim}\underset{n_2\rightarrow\infty}{\lim\sup }\ \underset{n_1\rightarrow\infty}{\lim\sup}\left|\left\{\left|m_{\varepsilon_{n_{1:3}}}-m\right|\geq\delta\right\}\right|=0\text{ on  }(0,T)\times\{v_0=0\}\text{ for all }\delta>0.\label{mconv0meas}
\end{align}
Finally, combining \eqref{mui}, \eqref{mconv0meas} and using the Vitali convergence theorem, we arrive at
\begin{align}
  \underset{n_3\rightarrow\infty}{\lim}\underset{n_2\rightarrow\infty}{\lim\sup}\ \underset{n_1\rightarrow\infty}{\lim\sup}\left\|m_{\varepsilon_{n_{1:3}}}-m\right\|_{L^1((0,T)\times\{v_0=0\})}=0.\label{cmstr0}
\end{align}
% \subsection*{Convergence for $m_{\varepsilon_{n_{1:3}}}$ (conclusion)}
% due to \eqref{cmstr}, \eqref{cmstr0}
% \begin{align}
%  \underset{n_3\rightarrow\infty}{\lim}\underset{n_2\rightarrow\infty}{\lim\sup}\underset{n_1\rightarrow\infty}{\lim\sup}\left\|m_{\varepsilon_{n_{1:3}}}-m\right\|_{L^1((0,T)\times\Omega)}=0;\label{cmL1}
% \end{align}
\subsection*{Convergence for $p_{\varepsilon_{n_{1:3}}}$ in \eqref{cpe}-\eqref{eq2e}}
We may consider \eqref{cpe}-\eqref{eq2e} together with the corresponding initial conditions as an abstract ODE system with respect to the variables $p_{\varepsilon_{n_{1:3}}}$ and $v_{\varepsilon_{n_{1:3}}}$ regarding $m_{\varepsilon_{n_{1:3}}}$ as a parameter function:
\begin{align}
 \frac{d}{dt}\left(\begin{array}{c}p_{\varepsilon_{n_{1:3}}}\\v_{\varepsilon_{n_{1:3}}}\end{array}\right)=G\left(\left(\begin{array}{c}p_{\varepsilon_{n_{1:3}}}\\v_{\varepsilon_{n_{1:3}}}\end{array}\right),m_{\varepsilon_{n_{1:3}}}\right)\text{ in }\LOne, \nonumber
\end{align}
where the function $G:([0,C_p]\times[0,1])\times\R^+_0\rightarrow\R^+_0$ is clearly globally Lipschitz. Here $C_p$ is an upper bound for the family $\{p_{\varepsilon}\}$, 
compare \eqref{cpb}. 
Using the standard abstract ODE theory in $L^1$, which states that the solutions depend continuously upon parameters and initial data, we conclude with \eqref{pk0}-\eqref{vk0} and \eqref{cmstr}, \eqref{cmstr0} that
\begin{align}
  &p_{\varepsilon_{n_{1:3}}}\underset{n_{1:3}\rightarrow\infty}{\rightarrow}p\text{ in  }L^1((0,T)\times\{v_0>0\}),\label{cpstr}\\
 & p_{\varepsilon_{n_{1:3}}}\underset{n_{1:3}\rightarrow\infty}{\rightarrow}p\text{ a.e. in }(0,T)\times\{v_0>0\},\label{pae}\\
 &\underset{n_3\rightarrow\infty}{\lim}\underset{n_2\rightarrow\infty}{\lim\sup}\ \underset{n_1\rightarrow\infty}{\lim\sup}\left\|p_{\varepsilon_{n_{1:3}}}-p\right\|_{L^1((0,T)\times\{v_0=0\})}=0,\label{cpstr0}
\end{align}
and $m,p,v$ solve the original equations \eqref{cp}-\eqref{v} and satisfy the initial conditions in $L^1$-sense, as stated in {\it Definition \ref{Defweak}}.
\subsection*{Convergence in {\eqref{weakfe}}}
% {\cb Hier stimmt die Referenz nicht, wir brauchen eine schwache Formulierung f\"ur die Approximation}
In order to finish the proof of {\it Theorem~\ref{maintheo}}, it remains to check that the triple $(m,p,v)$, which we obtained above by means of our limit procedure, satisfies the weak formulation \eqref{weakf}. { For this purpose, we need to pass to the limit  in the weak formulation \eqref{weakfe}. Taking $\varepsilon:=\varepsilon_{n_{1:3}}$, we have that
\begin{align}
   &\int_{\Omega}m_{\varepsilon_{2,n_2}0}\varphi\,dx\psi(0)-\int_0^T\int_{\Omega}m_{\varepsilon_{n_{1:3}}}\varphi\,dx\psi'\,dt\nonumber\\
   =&\int_0^T\int_{\Omega}-\varepsilon_{n_{1:3}}\nabla m_{\varepsilon_{n_{1:3}}}\cdot\nabla\varphi\psi-I_{\varepsilon_{n_{1:3}}}\cdot\nabla\varphi\psi+\kappa_v\nabla\left(\int_0^t\frac{v_{\varepsilon_{n_{1:3}}}m_{\varepsilon_{n_{1:3}}}}{1+v_{\varepsilon_{n_{1:3}}}}\,d\tau\right)\cdot\nabla\varphi\psi'\nonumber\\
   &\quad\quad\quad+\left(-\alpha m_{\varepsilon_{n_{1:3}}}+\beta v_{\varepsilon_{n_{1:3}}} p_{\varepsilon_{n_{1:3}}}\right)\varphi\psi\,dx\,dt,\label{weakfen}
  \end{align}
  where in order to shorten the notation we introduced 
  \begin{align}
   I_{\varepsilon_{n_{1:3}}}
   :=&\left(\frac{ \kappa_m c_{\varepsilon_{n_{1:3}}}}{1+v_{\varepsilon_{n_{1:3}}}c_{\varepsilon_{n_{1:3}}}}+\frac{\kappa_v}{1+v_{\varepsilon_{n_{1:3}}}}\right)2 v_{\varepsilon_{n_{1:3}}} ^{\frac{1}{2}}(m_{\varepsilon_{n_{1:3}}}+1)^{\frac{1}{2}}\nonumber\\
   &\cdot\left(\nabla\left(v_{\varepsilon_{n_{1:3}}}^{\frac{1}{2}}(m_{\varepsilon_{n_{1:3}}}+1)^{\frac{1}{2}}\right)-(m_{\varepsilon_{n_{1:3}}}+1)^{\frac{1}{2}}\nabla v_{\varepsilon_{n_{1:3}}}^{\frac{1}{2}}\right)\label{Ie}\\
   =&\frac{\kappa_m v_{\varepsilon_{n_{1:3}}}c_{\varepsilon_{n_{1:3}}}}{1+v_{\varepsilon_{n_{1:3}}}c_{\varepsilon_{n_{1:3}}}}\nabla m_{\varepsilon_{n_{1:3}}}+\frac{\kappa_v v_{\varepsilon_{n_{1:3}}}}{1+v_{\varepsilon_{n_{1:3}}}}\nabla m_{\varepsilon_{n_{1:3}}}.\label{Ie2}
  \end{align}
Observe that the representations \eqref{Ie} and \eqref{Ie2}
 coincide due to the chain and product rules. 
But for $I_{\varepsilon_{n_{1:3}}}$}, the convergence of the terms in \eqref{weakfen} can be obtained with standard tools using  \eqref{cpb}, \eqref{mui}, \eqref{reld}, \eqref{mk0}, \eqref{vae}, 
\eqref{vae_}, \eqref{cpstr}. We thus leave these details aside and concentrate on the weak $L^1$-limit for $I_{\varepsilon_{n_{1:3}}}$. To start with, \eqref{ddui}, \eqref{ta1ui} and {\eqref{Ie2}} imply that 
\begin{align}
\left\{I_{\varepsilon_{n_{1:3}}}\right\} \text{ is uniformly integrable in }(0,T)\times\Omega.\label{Ieui}                                                                                                                                                                                                                                                                                                                                                                                           \end{align}
Hence, the Dunford-Pettis theorem applies and yields the existence of such limit:
\begin{align}
I_{\varepsilon_{n_{1:3}}}\underset{n_{1:3}\rightarrow\infty}{\rightharpoonup}\tilde{I}\text{ in }L^1((0,T)\times\Omega).                 
\end{align}
We claim that $\tilde{I}$ can be obtained by simply dropping the index $\varepsilon_{n_{1:3}}$ everywhere in  \eqref{Ie}. We observe that \eqref{Ie} admits the following reformulation: 
\begin{align}
 I_{\varepsilon_{n_{1:3}}}=&I_1\left(m_{\varepsilon_{n_{1:3}}},p_{\varepsilon_{n_{1:3}}},v_{\varepsilon_{n_{1:3}}}\right) \nabla\left(I_3\left(m_{\varepsilon_{n_{1:3}}},v_{\varepsilon_{n_{1:3}}}\right)\right)+I_2\left(m_{\varepsilon_{n_{1:3}}},p_{\varepsilon_{n_{1:3}}},v_{\varepsilon_{n_{1:3}}}\right)\nabla v_{\varepsilon_{n_{1:3}}}^{\frac{1}{2}},\label{I123}
\end{align}
where $I_1,I_2:\R^3\mapsto\R$, $I_3:\R^2\mapsto\R$ are continuous functions. Since $I_1(\cdot,\cdot,0)\equiv I_2(\cdot,\cdot,0)\equiv0$, it holds with \eqref{vnposn}, \eqref{I123} that
\begin{align}
 I_{\varepsilon_{n_{1:3}}}=0\text{ a.e. in }(0,T)\times\{v_{\varepsilon_30}=0\}.\label{I0}
\end{align}
Combining \eqref{Ieui}, \eqref{I0} with   \eqref{vnposn}, \eqref{vnpos} and  property \eqref{star} and passing  to the limit in the measure on $(0,T)\times\{v_0=0\}$, we obtain that, as expected,
\begin{align}
 \tilde{I}=0\text{ a.e. in }(0,T)\times\{v_0=0\}.\nonumber
\end{align}
Further, we have due to \eqref{vae}, \eqref{mae}, \eqref{pae}, and the continuity of $I_1,I_2,I_3$ that
\begin{align}
 &I_1\left(m_{\varepsilon_{n_{1:3}}},p_{\varepsilon_{n_{1:3}}},v_{\varepsilon_{n_{1:3}}}\right)\underset{n_{1:3}\rightarrow\infty}{\rightarrow}I_1(m,p,v)\text{ a.e. in  }(0,T)\times\{v_0>0\},\label{I1}\\
 &I_2\left(m_{\varepsilon_{n_{1:3}}},p_{\varepsilon_{n_{1:3}}},v_{\varepsilon_{n_{1:3}}}\right)\underset{n_{1:3}\rightarrow\infty}{\rightarrow}I_2(m,p,v)\text{ a.e. in  }(0,T)\times\{v_0>0\},\label{I2}\\
 &I_3\left(m_{\varepsilon_{n_{1:3}}},v_{\varepsilon_{n_{1:3}}}\right)\underset{n_{1:3}\rightarrow\infty}{\rightarrow}I_3(m,v)\text{ a.e. in  }(0,T)\times\{v_0>0\}.\label{I3}
\end{align}
Using \eqref{d1ui}, \eqref{I3} and  the Vitali convergence theorem, we obtain that
\begin{align}
 &I_3\left(m_{\varepsilon_{n_{1:3}}},v_{\varepsilon_{n_{1:3}}}\right)\underset{n_{1:3}\rightarrow\infty}{\rightarrow}I_3(m,v)\text{ in  }L^1((0,T)\times\{v_0>0\}).\label{I3L1}
\end{align}
Together with \eqref{dui}, this yields by using the Dunford-Pettis theorem  that
\begin{align}
 \nabla\left(I_3\left(m_{\varepsilon_{n_{1:3}}},v_{\varepsilon_{n_{1:3}}}\right)\right)\underset{n_{1:3}\rightarrow\infty}{\rightharpoonup}\nabla(I_3(m,v)){ \text{ in  }L^1((0,T)\times\{v_0>0\}).}  \label{nI3}
\end{align}
Finally, combining \eqref{nv2}, \eqref{I1}, \eqref{I2}, \eqref{nI3} and using {\it Lemma \ref{LemA2}}, we arrive at
\begin{align}
 &I_1\left(m_{\varepsilon_{n_{1:3}}},p_{\varepsilon_{n_{1:3}}},v_{\varepsilon_{n_{1:3}}}\right) \nabla\left(I_3\left(m_{\varepsilon_{n_{1:3}}},v_{\varepsilon_{n_{1:3}}}\right)\right)+I_2\left(m_{\varepsilon_{n_{1:3}}},p_{\varepsilon_{n_{1:3}}},v_{\varepsilon_{n_{1:3}}}\right)\nabla v_{\varepsilon_{n_{1:3}}}^{\frac{1}{2}}\nonumber\\
 \underset{n_{1:3}\rightarrow\infty}{\rightharpoonup}
 &I_1\left(m,p,v\right) \nabla\left(I_3\left(m,v\right)\right)+I_2\left(m,p,v\right)\nabla v^{\frac{1}{2}}\text{ in }L^1((0,T)\times\{v_0>0\}).  \nonumber
\end{align}
The proof of {\it Theorem~\ref{maintheo}} is thus completed. 
% \begin{Remark}\label{Remv01}
%  If $||v_0||_{L^{\infty}(O)}<1$, then  the taxis flux term does belong to $L^1(O)$. {\cb Irgendwie kommt mir diese Bemerkung unfertig vor, was willst Du damit sagen?}
% \end{Remark}

\section{Numerical Simulations}\label{numerics}

We discretize the PDE-ODE-ODE system \eqref{haptoGG} using a local mass
conservative and monotone finite volume method. We use the
software package Dune \cite{dunegridpaperI:08,dunegridpaperII:08,ISTL,ISTLParallel} and consider on the domain
$\Omega = (0,1)^2$ the structured quadrilateral grid Yaspgrid therein.

\subsection{Implementation}

Let $\mathcal C$ be the set of computational cells in the grid and denote
by $\mathcal E(c)$ the (inner) edges of the grid surrounding a cell $c$.
Then we approximate the vector $u = (m,p,v)^T$ in the space $\mathcal P_0^3$,
so the restriction of $u$ on a computational cell $c$ is a constant vector.
Due to
the nonlinearity of the system it is favorable to employ IMEX-splitting
schemes, so we may handle one part of the system implicitly and another
part explicitly. The reaction part 
\begin{equation*}
  \partial_t \hat u = \begin{pmatrix}
    -\alpha \hat u_1 + \beta \hat u_2 \hat u_3\\
    \alpha \hat u_1 - \beta \hat u_2 \hat u_3 + \mu_p \hat u_2 ( 1 - (\hat u_1 + \hat u_2) - \eta \hat u_3)\\
    \mu_v \hat u_3 ( 1 - \hat u_3) - \lambda \hat u_3 \hat u_1
  \end{pmatrix} = f(\hat u)
\end{equation*}
of the system \eqref{haptoGG} is cell-wise a simple ODE, which we solve via an explicit 4th order
Runge-Kutta method.

\noindent
For the convection-diffusion part we have
\begin{equation}\label{diffusionPart}
  \partial_t \tilde u - \begin{pmatrix}
    \nabla \cdot \left(\frac{\kappa_m \tilde u_3 (\tilde u_1 + \tilde u_2)}{1 + \tilde u_3 (\tilde u_1 + \tilde u_2)} \nabla \tilde u_1 - \left(\frac{\kappa_v}{(1 + \tilde u_1 + \tilde u_2)^2} \nabla \tilde u_3\right) \tilde u_1 \right)\\
    0\\
    0\end{pmatrix} = 0.
\end{equation}
The discretization in space now takes place with the aid of two-point flux
approximations as in \cite{Eymard_et_al_FiniteVolumeMethods}.
First we define
the diffusion coefficient $D(u) = \frac{\kappa_m u_3 (u_1 + u_2)}{1 + u_3 ( u_1 + u_2)}$ and the drift velocity $V(u) = \frac{\kappa_v}{(1 + u_1 + u_2)^2}$.
The convection velocity $V(u) \nabla u_3$ and the diffusion term
$D(u) \nabla u_1$ have both the same structure, therefore we use the same space
discretization. Hence we will only present the diffusive flux discretization
in detail.
We may integrate \eqref{diffusionPart} over a computational cell
$c$ by employing the Gauß theorem
\begin{equation*}
  \partial_t u|_c = \sum\limits_{e \in \mathcal E(c)} F_c^e + V_c^e \ (u_1)_e^+,
\end{equation*}
where $F^e$ is the approximation of the diffusive flux  and
$V^e$ is an approximation of the drift velocity though an edge $e$.
The symbol $(u_1)_e^+$ stands for a
simple upwinding scheme \cite{Eymard_et_al_FiniteVolumeMethods}.
To get a locally
mass conservative method, we require that for each edge $e$ between cells
$c$ and $c'$ we have $F_c^e + F_{c'}^e = 0$, as well as $V_c^e + V_{c'}^e = 0$.
This gives the possibility to
resolve the edge variables and for an edge $e$ between $c$ and $c'$ we have
\begin{equation*}
  F_c^e = \frac{D(u)|_c D(u)|_{c'}}{D(u)|_c + D(u)|_{c'}} ((u_1)|_{c'} - (u_1)|_c) \frac{2 |e|}{d(c,c')}.
\end{equation*}
The drift velocity is computed in the same way. Now denote by $\mathcal F(u)$ the space discretized convective and diffusive flux terms and let the timestep of our scheme be $\Delta t$. Then we resolve the
reaction terms explicity (these are cell-wise ODEs) with a Runge-Kutta
method (denoted by its numerical flux $\Phi_{RK}$),
while the convection-diffusion part will be handled via an implicit
Euler step: %since explicit methods require hard restrictions on the timestep in order to end up with a stable method:
\begin{equation}\label{FullyDiscreteScheme}
  u^{k+1} + \Delta t \mathcal F(u^{k+1}) = u^k + \Phi_{RK}(u^k).
\end{equation}
We solve the previous equation \eqref{FullyDiscreteScheme} by the
classical Newton method.

\subsection{Results}

We have to select initial conditions. Therefore we assume a grate-like
initial condition for $v$ and define the following sets:
\begin{align*}
  S_1 &= \{x \in \mathbb R^2 | x_2 \in (0.35,0.45)\}\\
  S_2 &= \{x \in \mathbb R^2 | x_2 \in (0.7,0.8)\}\\
  S_3 &= \{x \in \mathbb R^2 | \ |x_1 - \hat x| < 0.01,\text{ for }\hat x \in \{0.4,0.45,0.5,0.55,0.6,0.65\}\}\\
  S_4 &= \{x \in \mathbb R^2 | \ |x_1 - x_2 -\hat x| < 0.01, \text{ for }\hat x \in \{-0.2, -0.1, 0.0\}\}\\
  S_5 &= \{x\in \mathbb R^2 | \ |x_1 - 0.5 \cdot x_2 - \hat x| < 0.01, \text{ for }\hat x \in \{0.5,0.6\}\}
\end{align*}
Then we select the intuitive initial value for the tissue fibers as
\begin{equation*}
  \tilde v_0 = 0.9 \cdot \mathbf 1_{\Big \{x \in \bigcup\limits_{i=1}^5 S_i\Big \}},
\end{equation*}
Now we need to think about the initial conditions for the tumor variables.
We observe that migrating tumor cells (variable $m$) will pass into the proliferating regime if no tissue is available (at least it is highly improbable to find a migrating cell 
in absence of tissue fibers). This is to be incorporated into the initial condition for
$m$. For the initial population of proliferating tumor cells, however,
we do not have the tissue dependence, so we may also select initial
conditions for $p$ in absence of tissue. Due to the fact, however, that
proliferating cells do not migrate (go-or-grow dichotomy), we have to assume a small compact
support. We use random perturbations of the initial conditions to simulate
the effect of non-homogeneous tumor cell distributions. With all these considerations we select the initial conditions
for the cell variables in the form
\begin{align*}
  m_0(x) &= \mathbf 1_{\{|x-x_0|^2 < 0.02\}} \cdot \min(0.5 \cdot \Psi_{0.05}(|x-x_0|^2 + d), 1.0) \cdot \mathbf 1_{\Big \{x \in \bigcup\limits_{i=1}^5 S_i\Big \}}\\
  p_0(x) &= \mathbf 1_{\{|x-x_0|^2 < 0.01\}} \cdot \min(0.8 \cdot \Psi_{0.1}(|x-x_0|^2 + d), 1.0),
\end{align*}
where $x_0 = \left(\frac{1}{2}, \frac{1}{2}\right)^T$ and
\begin{equation*}
  \Psi_\sigma(s) = \frac{1}{2\pi\sigma} \exp\left(-\frac{s}{2 \sigma^2}\right).
\end{equation*}
The symbol $d$ in the initial conditions stands for the random perturbation.
We used here a uniform $\mathcal U(-0.01,0.04)$ distribution.
We are not done in the initial values section, because due to the
dissolving of the tissue fibers caused by the migrating cells, we have to
modify the initial values for $v$ a bit:
\begin{equation*}
  v_0 = \max\left(\tilde v_0 - (m_0 + p_0),0.0\right).
\end{equation*}

The remaining task is to select the parameters involved in the
model. Some of them are available from literature, but for the diffusion coefficient $\kappa_m$ and the
haptotactic coefficient $\kappa_v$ we select higher values for the diffusion (compared to the previous papers
\cite{eks15,ZSU,SSU}), as the migratory behavior of the cells is diffusion dominated. 
The tissue is distributed in a quite inhomogeneous way, however on a tissue fiber (or fiber bundle) the
material is homogeneous, meaning that the tissue gradient $\nabla v$
and whence the haptotaxis is vanishing. Nevertheless, haptotaxis is not negligible, as it describes the guidance of cell migration by the tissue fibers (dissolved or not). The concrete parameter selection is summarized in Table \ref{tabelle}.

\begin{table}[ht]
\begin{center}
\begin{tabular}{|c| c|c||c|c|c| }
\hline
Parameter & Value & \text{Source} & \text{Parameter} & Value & \text{Source} \\
\hline
  $\alpha$ & 0.01 & \cite{eks15} &$\mu_p$ & 0.3 &\cite{ZSU,SSU} \\
  $\beta$ & 0.2 & \cite{eks15} &$\mu_v$ & 0.021 &\cite{ZSU,SSU} \\
  $\kappa_m$ & 0.1 & estimated &$\eta$ & 1.75& \cite{ZSU,SSU}\\
  $\kappa_v$ & 0.1 & estimated & $\lambda$ & 0.1& \cite{ZSU,SSU}\\
  \hline
\end{tabular}
\vspace*{1mm}
\caption{Parameters used in the model.}\label{tabelle}
\end{center} 
\end{table}
\noindent

\noindent
The grid we use is a triangulation of the unit cube in two dimensions,
with 200 cells in each direction. So we also have to select a small
time step $\Delta t$. In our calculations we used $\Delta t = 0.01$ and
simulated the equation up to time $1000$.

\begin{figure}[h!]
  \centering
  \caption{Simulation results}\label{figure-ueberhaupt}
  \begin{subfigure}[b]{\textwidth}
    \includegraphics[width = .99\textwidth]{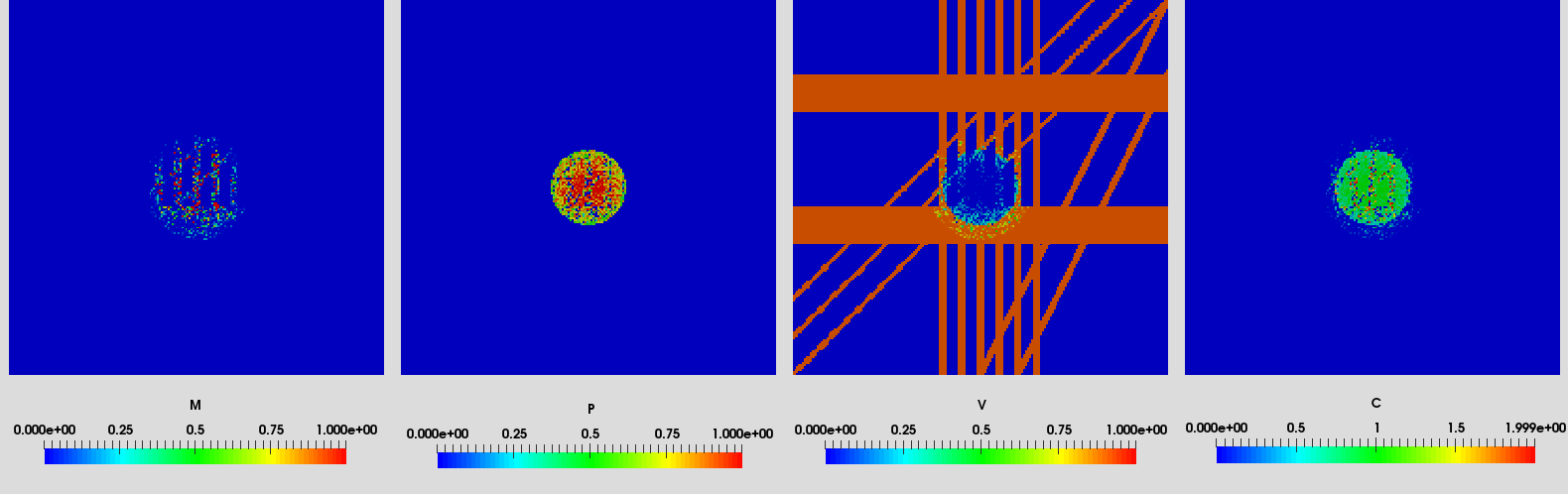}
    \caption{Initial condition. From left to right: migrating cells $m$, proliferating cells $p$, tissue $v$, overall tumor $c=m+p$.}
  \end{subfigure}\\
  \begin{subfigure}[b]{\textwidth}
    \includegraphics[width = .99\textwidth]{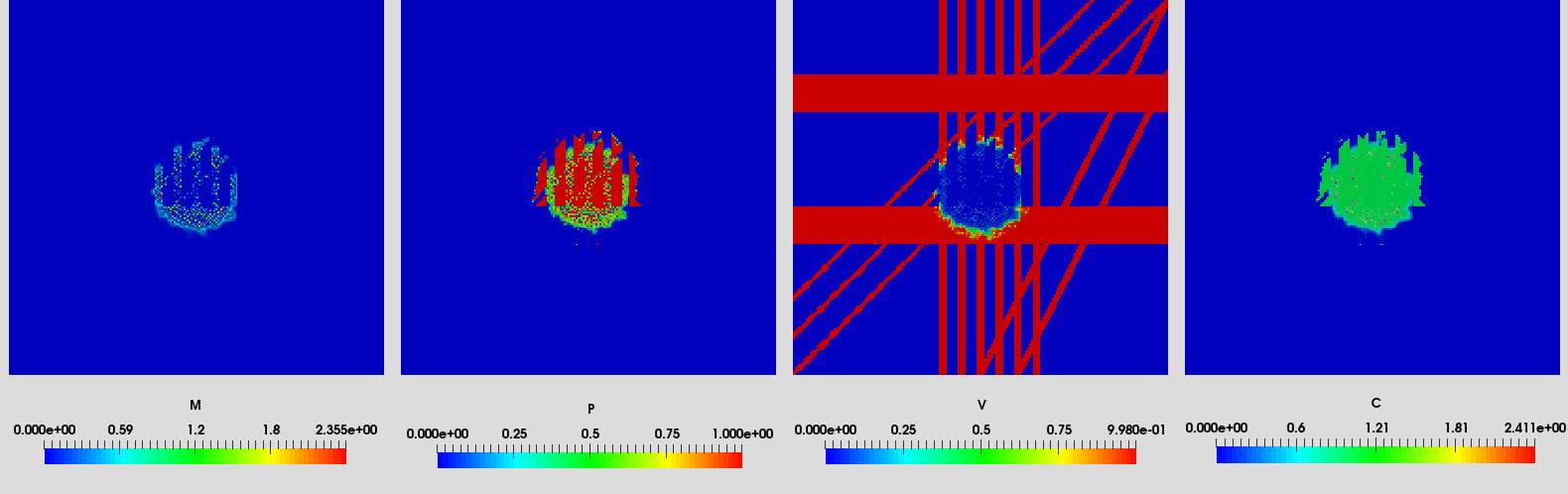}
    \caption{Simulation at time 200}
  \end{subfigure}\\
  \begin{subfigure}[b]{\textwidth}
    \includegraphics[width = .99\textwidth]{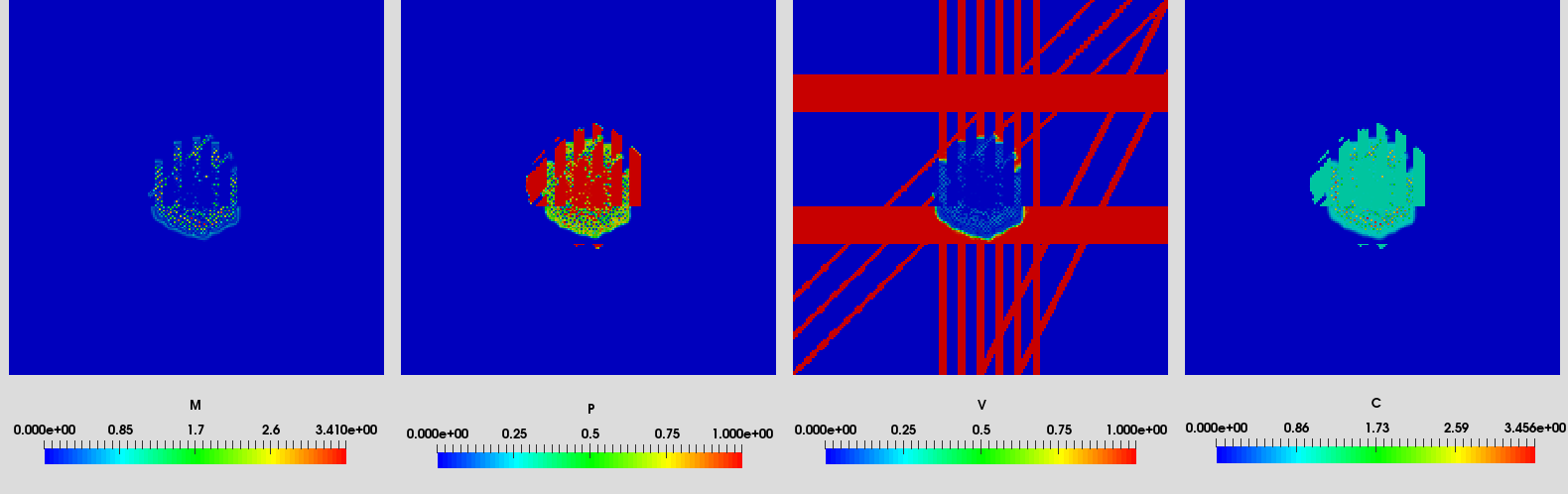}
    \caption{Simulation at time 400}
  \end{subfigure}
  \begin{subfigure}[b]{\textwidth}
    \includegraphics[width = .99\textwidth]{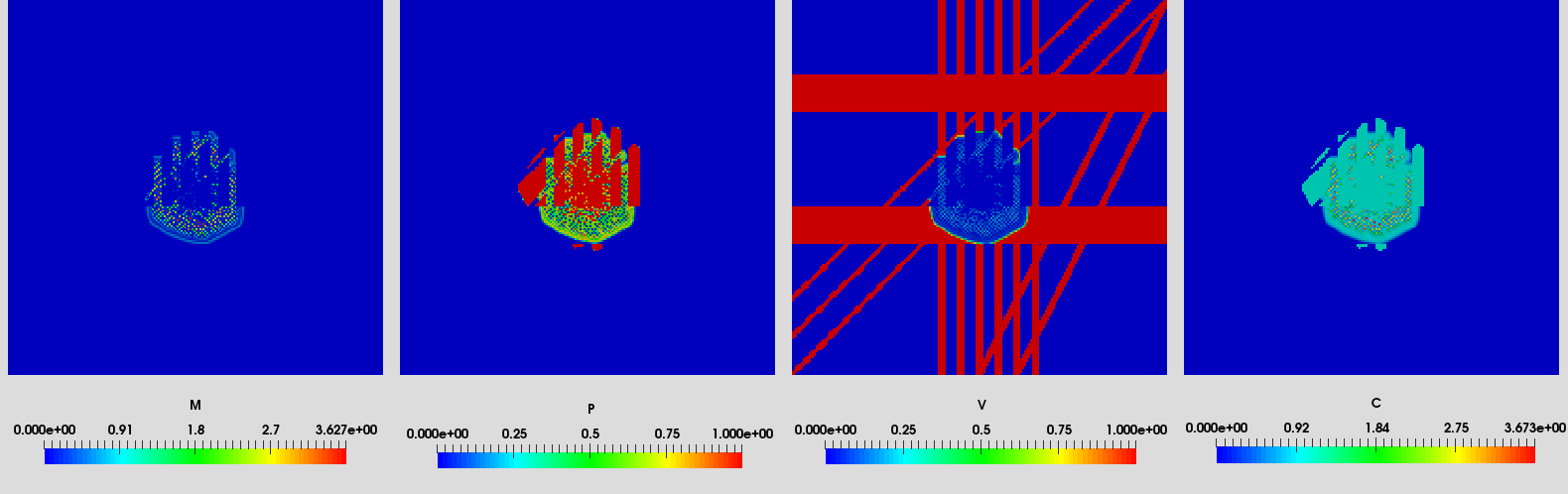}
    \caption{Simulation at time 600}
  \end{subfigure}
  %% \begin{subfigure}[b]{\textwidth}
  %%   \includegraphics[width = .99\textwidth]{OtherViewTime800}
  %%   \caption{Simulation to time 800}
  %% \end{subfigure}
\end{figure}

\noindent
{\it Figure \ref{figure-ueberhaupt}} shows the simulation results. The comparison between the evolution of migrating and 
proliferating cells elicits the expected behavior: the migrating cells are predominant in the regions with high tissue density 
(it can be actually seen how they follow the tissue fibers -and degrade them), while the proliferating cells occupy the regions 
with very low tissue density. This is in agreement with the go-or-grow dichotomy and the degeneracy of the diffusion coefficient 
in equation \eqref{cm}: For $v=0$ (no tissue) the migrating cells stop and become proliferating cells. Moreover, the 
model is able to reproduce the often irregular shape of a tumor and the associated spread of cancer cells exhibiting  
various infiltrative (INF) patterns. According to the Japanese gastric association group \cite{japanese}, the latter provide a 
way to classify local invasiveness and tumor malignancy. In particular, {\it Figure \ref{figure-ueberhaupt}} exhibits some small 
'islands' of cell aggregates, transiently isolated from the main tumor, which then grow and merge again with the neoplastic 
cell mass. That tumor cells have an infiltrative spread, form fingering patterns, and closely follow the specific tissue 
structure has been recognized for 
many types of cancer; perhaps the most prominent example featuring these characteristics are gliomas, see e.g. 
\cite{coons,gerstneretal2010,giese-etal96,giesewestphal1996}. This behavior has also been confirmed by several 
models in a different mathematical framework, but still relying on the go-or-grow dichotomy and leading to related 
reaction-diffusion-taxis equations \cite{eks15,hs15}. Like those models, the present setting allows to account for 
tumor heterogeneity w.r.t. the migratory/proliferative phenotypes of the constituent cells. As mentioned in the Introduction,
this heterogeneity also reflects in the differentiated therapeutic response, an essential issue in therapy 
planning and assessment. Including therapy effects like e.g., in \cite{SSU} can be easily addressed in this context as well. While  
current biomedical imaging only allows to determine the gross tumor volume, such models open the way to 
 provide an (although 
imperfect) estimation of the tumor composition upon relying on the patient-specific tissue architecture and to correspondingly 
predict the extent of the neoplastic tissue.

\noindent
Another interesting observation is that the amount of migrating cells increases with advancing time. This might suggest a 
possible blow-up; we recall that this issue remains open from an analytical point of view.

%\section{Discussion}

% ...
% \begin{align}
%  |\{m_{\varepsilon}=0\}\backslash(\{p_{\varepsilon}=0\}\cup\{v_{\varepsilon_30}=0\})|=0
% \end{align}
% \phantomsection
% \printbibliography
\begin{appendices}
\section{}\label{}
In this section we collect several auxiliary results on member-by-member products used above. We begin with a lemma which deals with the uniform integrability of member-by-member products.
\begin{Lemma}[Uniform integrability for  products]\label{UiL2}
 Let $\Omega$ be a measurable subset of $\R^N$ with finite measure and $I$ be some set. Let $\{ f_i\}_{i\in I}, \{ g_i\}_{i\in I}\subset L^2(\Omega)$ be two families such that $\{ |f_i|^2\}_{i\in I}$ is  uniformly integrable and $\{ g_i\}_{i\in I}$ is uniformly  bounded in $L^2(\Omega)$. Then the family $\{ f_ig_i\}_{i\in I}$ of member-by-member products is uniformly integrable.
\end{Lemma}
% \begin{proof}
% Indeed, due to the assumptions and the Cauchy-Schwartz inequality, it holds that
%  \begin{align}
% \underset{n\rightarrow\infty}{\lim}\underset{A\subset\Omega, |A|\leq\frac{1}{n}}{\sup}\underset{i\in I}{\sup}\left|\int_{A}f_ig_i\,dx\right|\leq& \underset{i\in I}{\sup}||g_i||\left(\underset{n\rightarrow\infty}{\lim}\underset{A\subset\Omega, |A|\leq\frac{1}{n}}{\sup}\underset{i\in I}{\sup}\int_{A}|f_i|^2\,dx\right)^{\frac{1}{2}}=0.\nonumber
% \end{align}
% Thus, the family of member-by-member products is uniformly integrable by definition.
% \end{proof}
\noindent
This well-known property can be readily proved by using the definition of the uniform integrability. We leave the details to the reader. 
\noindent
The following lemma is a generalization of the Lions lemma \cite[Lemma 1.3]{Lions} and the known result on weak-strong convergence for member-by-member products. %We recall its proof for the convenience of the reader.
\begin{Lemma}[Weak-a.e. convergence, \cite{ZSU}]\label{LemA1}
 Let $\Omega$ be a measurable subset of $\R^N$ with finite measure. Let $f,f_n:\Omega\rightarrow\R$, $n\in\N$ be measurable functions and $g,g_n\in L^1(\Omega)$, $n\in\N$. Assume further that $f_n\underset{n\rightarrow\infty}{\rightarrow} f$ a.e. in $\Omega$ and $g_n\underset{n\rightarrow\infty}{\rightharpoonup}g$, $f_ng_n\underset{n\rightarrow\infty}{\rightharpoonup}\xi$ in $L^1(\Omega)$. Then, it holds that $\xi=fg$ a.e. in $\Omega$.
 \end{Lemma}
\noindent
As was observed in \cite{ZSU}, a similar result  holds for sums of member-by-member products:
\begin{Lemma}[Weak-a.e. convergence for sums, \cite{ZSU}]\label{LemA2}
 Let $\Omega$ be a measurable subset of $\R^N$ with finite measure and let $L\in\N$. Let $f^l,f^l_n:\Omega\rightarrow\R$, $n\in\N$, $l\in\{1,...,L\}$, be measurable functions and $g^l,g^l_n\in L^1(\Omega)$, $n\in\N$, $l\in\{1,...,L\}$. Assume further that $f^l_n\underset{n\rightarrow\infty}{\rightarrow} f^l$ a.e. in $\Omega$ and $g^l_n\underset{n\rightarrow\infty}{\rightharpoonup}g^l$, $\sum_{l=1}^{L}f^l_ng^l_n\underset{n\rightarrow\infty}{\rightharpoonup}\xi$ in $L^1(\Omega)$. Then, it holds that $\xi=\sum_{l=1}^{L}f^lg^l$ a.e. in $\Omega$.
 \end{Lemma}
 \begin{Remark}
  Observe that, in {\it Lemma~\ref{LemA2}}, it is not required that the  sequences $\left\{f^l_ng^l_n\right\}_{n\in\N}$ themselves  are convergent for $l\in\{1,...,L\}$, but only their sum $\left\{\sum_{l=1}^{L}f^l_ng^l_n\right\}_{n\in\N}$.  Thus, the result is applicable in the cases where the convergence of individual sequences is either false or unknown.
 \end{Remark}
%  \noindent
% The proof of {\it Lemma~\ref{LemA2}} is very similar to the proof of {\it Lemma~\ref{LemA1}}. One only has to choose the sets $\Omega_k$ and  $\Omega_{k,m}$  independent of $l\in\{1,...,L\}$. We leave the remaining details to the reader.
\end{appendices}
\phantomsection
\printbibliography
\end{document}